\pdfoutput=1

\documentclass[
    english, 
    a4paper,
    12pt,
    oneside, 
    final,
    numbers=noenddot, 
    toc=bibliography,
]{scrbook}

\usepackage[utf8]{inputenc}
\usepackage{babel}

\usepackage{scrpage2}
\clearscrheadfoot

\ihead{\leftmark}
\ofoot{}
\setheadsepline{.4pt}
\ohead[]{\pagemark}

\setkomafont{pageheadfoot}{\normalfont}

\usepackage{titlesec}

\usepackage{microtype}

\usepackage{xspace}
\usepackage{rotating}

\usepackage{xcolor}

\newcommand{\inred}[1]{\xspace{\color{red} #1}\xspace} 

\usepackage[hyphens]{url}

\usepackage[parfill]{parskip}

\usepackage{amssymb, amsmath, bm}
\usepackage{amsthm}

\usepackage[expproduct=times]{siunitx} 

\usepackage{listings}
\usepackage{algorithm}
\usepackage{algpseudocode}

\usepackage{pgfplots} 
\newlength\fheight 
\newlength\fwidth 
\usepackage{tikz}

\usetikzlibrary{external,arrows, positioning, shapes,fit}

\usepackage{caption}
\usepackage{subcaption}

\usepackage{import}

\newlength{\myheight}

\usepackage{array}
\newcolumntype{L}[1]{>{\raggedright\let\newline\\\arraybackslash\hspace{0pt}}m{#1}}
\newcolumntype{C}[1]{>{\centering\let\newline\\\arraybackslash\hspace{0pt}}m{#1}}
\newcolumntype{R}[1]{>{\raggedleft\let\newline\\\arraybackslash\hspace{0pt}}m{#1}}

\usepackage{tabularx, etoolbox, booktabs}
\usepackage{longtable}
\usepackage{multirow} 
\newcommand{\rarray}[1]{\renewcommand{\arraystretch}{#1}} 
\usepackage{parcolumns}
\usepackage{paracol}
\usepackage{lipsum}

\usepackage{csquotes} 
\usepackage[backend=bibtex, style=alphabetic, firstinits=true]{biblatex}
\bibliography{literature}

\usepackage[acronym, toc, nonumberlist, nopostdot, nogroupskip]{glossaries}
\newacronym{ADSL}{ADSL}{Asymmetric Digital Subscriber Line}
\newacronym{FFT}{FFT}{fast fourier transformation}
\newacronym{dft}{DFT}{discrete fourier transformation}
\newacronym{idft}{IDFT}{inverse discrete fourier transformation}
\newacronym{rdft}{RDFT}{discrete fourier transformation for real input}
\newacronym{irdft}{IRDFT}{inverse discrete fourier transformation for real input}
\newacronym{NN}{NN}{neural network}
\newacronym{PSD}{PSD}{power spectral density}
\newacronym{blstm}{BLSTM}{bi-directional long short-term memory}
\newacronym{pca}{PCA}{principle component analysis}
\newacronym{mvdr}{MVDR}{minimum variance distortionless response}
\newacronym{gev}{GEV}{generalized eigenvector}
\newacronym{asr}{ASR}{automatic speech recognition}
\newacronym{ad}{AD}{algorithmic differentiation}
\newacronym{DNN}{DNN}{deep neural network}
\newacronym{ban}{BAN}{blind analytic normalization}
\newacronym{pesq}{PESQ}{perceptual evaluation of speech quality \cite{Rix2001PESQ}}
\newacronym{stft}{STFT}{short-time Fourier transform}
\newacronym{ff}{FF}{feedforward}
\newacronym{ibm}{IBM}{ideal binary mask}
\newacronym{bce}{BCE}{binary cross entropy}
\newacronym{wer}{WER}{word error rate}
\newacronym{snr}{SNR}{signal-to-noise ratio}
\newacronym{ml}{ML}{maximum likelihood}

\makeglossaries
\glsenablehyper

\usepackage[colorlinks=false, pdfborder={0 0 0}, linktoc=all,hidelinks]{hyperref}
\usepackage[capitalize]{cleveref}

\usepackage{soul}
\usepackage{calc}
\usepackage{placeins} 

\usepackage{enumitem}

\DeclareMathOperator{\argmaxOp}{argmax}
\newcommand\argmax[2]{\underset{#1}{\argmaxOp} \left\{ #2 \right\} }

\DeclareMathOperator{\E}{\mathbb{E}}
\DeclareMathOperator{\e}{e}

\makeatletter
\newcommand{\const}{\mathrm{const}\@ifnextchar.{}{.}}
\makeatother

\renewcommand{\vec}[1]{\ensuremath{\mathbf{#1}}}
\newcommand{\vecG}[1]{\ensuremath{\boldsymbol{#1}}}

\newcommand{\minusOne}{-1}
\newcommand{\minusH}{-\mathrm{H}}

\DeclareMathOperator{\hadamard}{\circ}  

\newcommand{\inv}{^{\minusOne}}
\newcommand{\transposed}{^{\mathrm{T}}}
\newcommand{\transposedinv}{^{-\mathrm{T}}}
\newcommand{\hermite}{^{\mathrm{H}}}
\newcommand{\hermiteinv}{^{\minusH}}

\newcommand*\der[2]{\frac{\partial #1}{\partial #2}}
\newcommand{\grad}[1]{\frac{\partial J}{\partial #1\conj}}
\newcommand{\Tr}{\mathrm{Tr}}
\newcommand{\conj}{^*}
\newcommand{\trace}[1]{\mathrm{Tr}\left(#1\right)}
\newcommand{\Complex}{\mathbb{C}}
\newcommand{\Real}{\mathbb{R}}
\renewcommand{\j}{\mathrm{j}}

\newcommand{\phixx}{{\vecG{\Phi}_{XX}}}
\newcommand{\phixxtilde}{{\vecG{\Phi}_{\tilde{X}\tilde{X}}}}
\newcommand{\phinn}{{\vecG{\Phi}_{NN}}}
\newcommand{\w}{{\vec{w}}}
\newcommand{\sumD}{\sum_{d=1}^{D}}
\newcommand{\sumT}{\sum_{t=1}^{T}}
\newcommand{\sumTc}{\sum_{\check{t}=1}^{T}}

\newcommand{\X}{{\vec{X}}}
\newcommand{\N}{{\vec{N}}}
\newcommand{\Y}{{\vec{Y}}}
\newcommand{\Ytilde}{{\vec{\tilde{Y}}}}

\newcommand{\Ltwo}[1]{\left|\left|#1\right|\right|}
\newcommand{\PCA}{\ensuremath{\mathrm{PCA}}}
\newcommand{\GEV}{\ensuremath{\mathrm{GEV}}}
\newcommand{\MVDR}{\ensuremath{\mathrm{MVDR}}}
\newcommand{\wrt}{w.r.t.\xspace}
\newcommand{\rv}{real-valued\xspace}
\newcommand{\cv}{complex-valued\xspace}

\newcommand{\nuVec}{\vecG{\nu}}
\newcommand{\sumF}{\sum_{f=1}^{F}}
\newcommand{\Norm}{\mathrm{Norm}}

\newcommand*\derz[1]{\frac{\partial #1}{\partial \zeta}}

\newcommand*{\principalComponent}[1]{\mathcal{P}\left(#1\right)}

\newcommand*{\llblacktriangle}[1][]{\tikz[x=0.7em, y=0.7em]\fill[#1] (0,0) -- (0,1) -- (1,0) -- cycle;\xspace}
\newcommand*{\urblacktriangle}[1][]{\tikz[x=0.7em, y=0.7em]\fill[#1] (1,1) -- (0,1) -- (1,0) -- cycle;\xspace}
\newcommand*{\lltriangle}[1][]{\tikz[x=0.7em, y=0.7em]\draw[#1] (0,0) -- (0,1) -- (1,0) -- cycle;\xspace}
\newcommand*{\urtriangle}[1][]{\tikz[x=0.7em, y=0.7em]\draw[#1] (1,1) -- (0,1) -- (1,0) -- cycle;\xspace}

\newcommand{\dft}[2][n\rightarrow f]{\underset{#1}{\mathrm{DFT}}\left\{#2\right\}}
\newcommand{\idft}[2][f\rightarrow n]{\underset{#1}{\mathrm{IDFT}}\left\{#2\right\}}		
\newcommand{\rdft}[2][n\rightarrow f]{\underset{#1}{\mathrm{RDFT}}\left\{#2\right\}}
\newcommand{\irdft}[2][f\rightarrow n]{\underset{#1}{\mathrm{IRDFT}}\left\{#2\right\}}

\globalcounter{equation}

\begin{document}
    
    \frontmatter

\begin{titlepage}
\begin{center}
  \begin{figure}\vspace{-1cm}
      \includegraphics[scale=0.45]{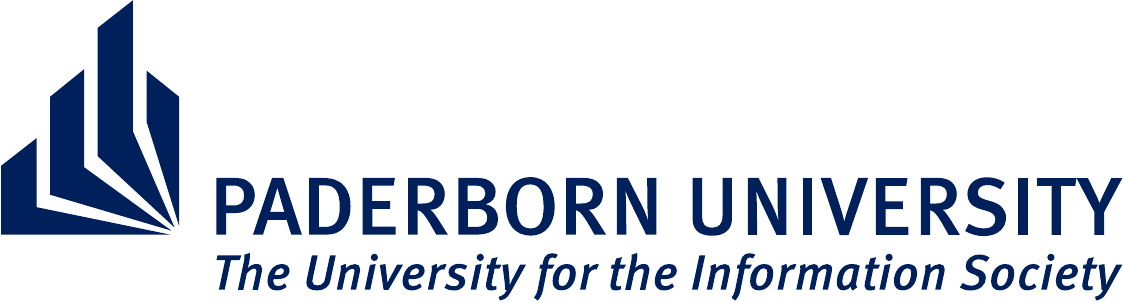} \hspace{4cm}
      \includegraphics[scale=1.2]{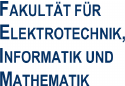}
  \end{figure}
$\:$ \newline
{\textbf{Institute for Electrical Engineering and Information Technology}}\\
Paderborn University\\
Department of Communications Engineering\\
Prof. Dr.-Ing. Reinhold Haeb-Umbach\\
\vfill
$\:$ \newline
\Large{\textbf{\textsf{Technical Report \vspace{1cm} \\ On the Computation of Complex-valued Gradients \\ with Application to \\Statistically Optimum Beamforming}}}
\vfill
\large{by\\}
\vfill
\Large{\large Christoph Boeddeker$^\star$, Patrick Hanebrink$^\star$, \\Lukas Drude, Jahn Heymann, Reinhold Haeb-Umbach\\
\texttt{\normalsize\{cbj,phane\}@mail.upb.de, \{drude,heymann,haeb\}@nt.upb.de}
}
\vfill
\renewcommand{\baselinestretch}{1}\normalsize
\end{center}
\vspace{0.5cm}
\footnoterule
{\footnotesize \hspace{\footnotesep} $^\star$Both authors contributed equally.}
\end{titlepage}

    \titleformat{\chapter}
    [display] 
    {\bfseries\huge} 
    {} 
    {0ex}
    {\vspace{-5ex}\titlerule\vspace{1.5ex}\filright~}
    [\vspace{1ex}\titlerule]
    

\chapter{Abstract}
\begin{center}
\begin{minipage}{0.6\textwidth}
This report describes the computation of gradients by algorithmic differentiation for statistically optimum beamforming operations.
Especially the derivation of complex-valued functions is a key component of this approach.
Therefore the real-valued algorithmic differentiation is extended via the complex-valued chain rule.
In addition to the basic mathematic operations the derivative of the eigenvalue problem with complex-valued eigenvectors is one of the key results of this report.
The potential of this approach is shown with experimental results on the CHiME-3 challenge database.
There, the beamforming task is used as a front-end for an ASR system. 
With the developed derivatives a joint optimization of a speech enhancement and speech recognition system \wrt the recognition optimization criterion is possible.
\end{minipage}
\end{center}
    
    \tableofcontents
    
    \mainmatter
    \titleformat{\chapter}
    [display] 
    {\bfseries\huge} 
    {} 
    {0ex}
    {\vspace{-5ex}\titlerule\vspace{1.5ex}\filright \thechapter~}
    [\vspace{1ex}\titlerule]

\chapter{Introduction}
\label{cha:introduction}

Over the years \glspl{NN} have become state-of-the-art in many speech signal processing, e.g. \gls{asr}.
Recently they also gained impact in speech enhancement.
There a \gls{NN} estimates masks which are used to calculate the \gls{PSD} matrices of the speech and noise signal \cite{merl2016mvdr,Heymann2015,Heymann2016}.
With these \gls{PSD} matrices the beamforming coefficients are calculated with a statistical model based approach.
Used with a strong \gls{asr} backend this system can improve the performance notably.

Looking at the complete system containing both speech enhancement and \gls{asr}, it is suboptimal to have two separate \glspl{NN} involved.
Both are trained with different objective functions which does not ensure, that the mask estimation is optimal in terms of speech recognition.
To overcome this problem the gradients of the model-based beamforming operations have to be known to use the backpropagation algorithm for optimizing \glspl{NN}.
Since the beamformer operates with the \gls{stft} of the noisy signal \cv operations are involved.
This requires \cv gradients of many operations, including the eigendecomposition.

The objective function used in this report is the output \gls{snr} of the beamformer, which is closer related to speech enhancement than the mask purity criterion.
The gradients are calculated using \gls{ad}, where a function is decomposed into a chain of elementary functions of which the gradient \wrt the desired quantity can be calculated.
In especially we extend many well-known \rv gradients to their \cv case.
Therefore different concepts of deriving those gradients are presented.
Then some examples are given, including the gradient of the \gls{dft}, the matrix inverse and the \cv eigenvalue problem.
Also the gradient of the Cholesky decomposition is given.

In \cref{cha:status} the system consisting of a beamformer and \gls{asr} part is shown.
Also the statistically beamforming operations are described.
The complex-valued backpropagation, gradient descent and the derivation of the backward gradient are explained in \cref{cha:complexbackprop}.
\cref{cha:blocks} shows the derivation of the backward gradients for many commonly used complex-valued operations. An overview of all gradients can be found in \cref{sec:blockSummary}.
The evaluation of the system with the derived gradients is done in \cref{cha:evaluation}.
\cref{cha:summary} summarizes the results and gives a short outlook.


\chapter{Model} 
\label{cha:status}

\section{Description of the system}
\label{sec:model}

The system contains a $D$-dimensional microphone array, that receives a noisy speech signal
\begin{align}
    \Y_{f, t} &= \X_{f, t} + \N_{f, t},
\end{align}
where the D channels are combined to an observation vector.
$\Y_{f, t}$, $\X_{f, t}$ and $\N_{f, t}$ are the \gls{stft} of noisy speech, clean speech and noise at frequency bin $f$ and time frame $t$, respectively.
The system is designed to extract the clean speech signal $X_{f,t,1}$ observed by the first microphone.
Here acoustic beamforming with the beamformer coefficient vector $\vec{w}_f$, which is assumed to be time-invariant, is used \cite[Ch. 10]{LiDeHaGo2015ASR}:
\begin{align}
    \hat{X}_{f,t,1} = \w_f^{\mathrm{H}} \cdot \Y_{f,t}. \label{eq:beamformer_apply}
\end{align}

To obtain the beamformer coefficient vector the speech and noise \gls{PSD} matrices have to be estimated.
This is done via a \gls{NN}, which is applied to each microphone channel separately.
It estimates a speech mask $M^{(X)}_{f, t, d}$ and a noise mask $M^{(N)}_{f, t, d}$ for each channel.
After combining (i.e. averaging) the $D$ speech and noise masks to a single speech and a single noise mask and assuming a zero-mean observation signal $\Y_{f, t}$ the \gls{PSD} matrices can be estimated with a \gls{ml} estimator \cite[Eq. 9.19]{Bishop2006}
\begin{align}
	{\vecG{\Phi}_{\nu\nu}}_f &= \frac{\sumT\left(\sumD M^{(\nu)}_{f, t, d}\right)\Y_{f, t} \cdot \Y_{f, t}\hermite}{\sumT \left(\sumD M^{(\nu)}_{f, t, d}\right)}, \label{eq:psd_calculation}
\end{align}
where $\nu \in \{X,N\}$ and where $(\cdot)\hermite$ denotes the conjugate transpose.
With these estimated \gls{PSD} matrices different statistically optimal beamformers can be computed (see \cref{sec:beamformers}).

An objective function calculated from the beamformer output is used to update the parameters of the \gls{NN}. In this report the negative \gls{snr} of the normalized beamformer output signal
\begin{align}
	J &= - 10 \cdot \log_{10}{\frac{P^{(X)}}{P^{(N)}}}, \label{eq:lossFkn}
\end{align}
with the power of the beamformed signal
\begin{align}
	P^{(\nu)} &= \frac{1}{T} \sumT \sumF |\w_{f}\hermite \cdot \nuVec^{(\Norm)}_{f, t}|^2, \label{eq:snr:loss}
\end{align}
applied to the normalized clean speech and noise signal
\begin{align}
	\nuVec^{(\Norm)}_{f, t} &= \frac{\nuVec_{f, t}}{\sqrt{\sumTc \Ltwo{\nuVec_{f, \check{t}}}^2}},
\end{align}
is used as the objective function. $T$ and $F$ are the total number of frames and frequency bins, respectively.
Every frequency contributes equally strong to the objective due to the normalization.

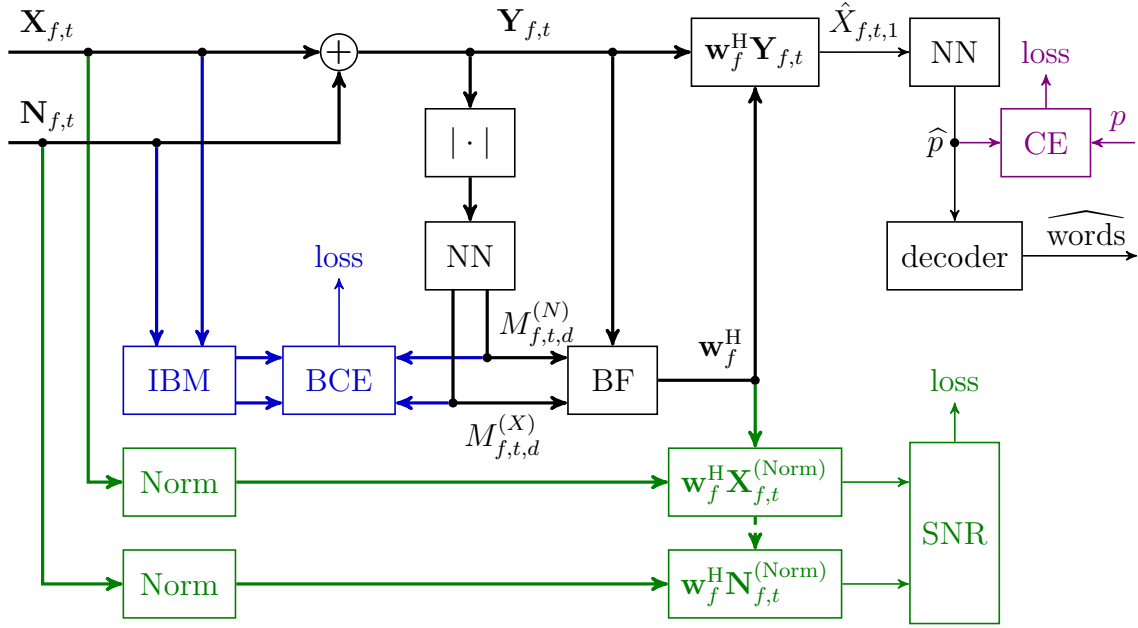
\begin{figure}
	\centering
	\definecolor{oldCostCol}{RGB}{0,0,200}
	\definecolor{newCostCol}{RGB}{0,128,0}
	\definecolor{werCostCol}{RGB}{128,0,128}

\begin{tikzpicture}[scale=1.5,
		>=stealth',semithick,auto,node distance = 3cm,
		block/.style={rectangle,draw,text centered, text width=0.9cm, minimum height=0.9cm},
		blockCostNew/.style={rectangle,draw,text centered, text width=1.2cm, minimum height=0.9cm, draw=newCostCol, text=newCostCol},
		blockCostOld/.style={rectangle,draw,text centered, text width=1.2cm, minimum height=0.9cm, draw=oldCostCol, text=oldCostCol},
		blockCostWer/.style={rectangle,draw,text centered, text width=0.9cm, minimum height=0.9cm, draw=werCostCol, text=werCostCol},
		branch/.style = {circle,inner sep=0pt,minimum size=1mm,fill=black,draw=black},
		sum/.style = {draw, circle, inner sep = 0pt},
		clear/.style = {inner sep = 0pt}
	]

	\newcommand{\beamformOp}{$\w_f \hermite (\cdot)$}
	
	\pgfmathsetmacro{\cOne}{0} %
	\pgfmathsetmacro{\cTwo}{\cOne+1.5} 
	\pgfmathsetmacro{\cThree}{\cTwo+1.4} 
	\pgfmathsetmacro{\cFour}{\cThree+1.15} 
	\pgfmathsetmacro{\cFife}{\cFour+1.25} 
	\pgfmathsetmacro{\cSix}{\cFife+1.25} 
	\pgfmathsetmacro{\cSeven}{\cSix+1.75} 
	\pgfmathsetmacro{\cEight}{\cSeven+0.8} 
	\pgfmathsetmacro{\cNine}{\cEight+0.8} 
	
	\pgfmathsetmacro{\rOne}{0} %
	\pgfmathsetmacro{\rTwo}{\rOne+0.9} 
	\pgfmathsetmacro{\rThree}{\rTwo+0.9} 
	\pgfmathsetmacro{\rFour}{\rThree+1.1} 
	\pgfmathsetmacro{\rFife}{\rFour+1} 
	\pgfmathsetmacro{\rSix}{\rFife+0.8} 


	\coordinate (x) at (\cOne,\rSix);
	\node[above right] at (x){$\X_{f,t}$};
	\coordinate (n) at (\cOne,\rFife);
	\node[above right] at (n){$\N_{f,t}$};
	\node[sum](sum) at (\cThree,\rSix){\textbf{+}};
	\node[branch](br5)at (\cFour,\rSix){};
	\node[block](abs) at (\cFour,\rFife){$|\cdot|$};
	\node[block](NN1)at (\cFour,\rFour){NN};
	\node[block](bf) at (\cFife,\rThree){BF};
	\node[branch](br6)at (\cFife,\rSix){};
	\node[block,text width=1.4cm](aplBf1) at (\cSix,\rSix){$\w_f \hermite \Y_{f,t}$};
	\node[block](NN2)at (\cSeven,\rSix){NN};
	\node[block, text width=1.5cm](decode) at(\cSeven,\rFour){decoder};
	\coordinate (wordHat) at (\cNine,\rFour);
	\node[above left]at (wordHat){$\widehat{\text{words}}$};

	\node[branch](br3)at (\cTwo-0.2,\rFife){};
	\node[branch](br4)at (\cTwo+0.2,\rSix){};
	
	\node[blockCostOld](ibm) at(\cTwo,\rThree){IBM};
	
	\node[blockCostOld](bce) at (\cThree,\rThree){BCE};
	\node[text=oldCostCol](lossOld) at (\cThree,\rFour){loss};
	\node[branch](br8) at (\cFour+0.15,\rThree+0.2){};
	\node[above right]at (\cFour+0.15,\rThree+0.2) {$M_{f,t,d}^{(N)}$};
	\node[branch](br9) at (\cFour-0.15,\rThree-0.2){};
	\node[below right] at (\cFour-0.15,\rThree-0.2) {$M_{f,t,d}^{(X)}$};

	\node[branch](br1)at (\cTwo-1.2,\rFife){};
	\node[branch](br2)at (\cTwo-0.8,\rSix){};
	\node[blockCostNew](norm1) at(\cTwo,\rTwo){Norm};
	\node[blockCostNew](norm2) at(\cTwo,\rOne){Norm};
	\node[branch](br7) at (\cSix,\rThree){};
	\node[blockCostNew, text width=2.0cm](aplBf2) at (\cSix,\rTwo){$\w_f \hermite \X_{f,t}^{(\mathrm{Norm})}$};
	\node[blockCostNew, text width=2.0cm](aplBf3) at (\cSix,\rOne){$\w_f \hermite \N_{f,t}^{(\mathrm{Norm})}$};
	\node[blockCostNew, text width = 0.9cm, minimum height = 2.4cm](snr) at (\cSeven,\rOne/2+\rTwo/2){SNR};
	\node[text=newCostCol](lossNew) at (\cSeven,\rThree){loss};

	\node[blockCostWer](ce) at (\cEight,\rFife){CE};
	\node[branch](br10) at (\cSeven,\rFife){};
	\node[left]at (\cSeven,\rFife){$\widehat{p}$};
	\node[clear](p) at (\cNine,\rFife){};
	\node[above left,text=werCostCol] at (p){$p$};
	\node[text=werCostCol](lossWer) at (\cEight,\rSix){loss};

	\path [{}-{>}]
		(x)  edge[very thick] (sum)
		(sum) edge[very thick] node{$\Y_{f,t}$} (aplBf1)
		(decode) edge (wordHat)
		(NN2) edge (decode)
		(br5) edge[very thick] (abs)
		(abs) edge[very thick] (NN1)
		(br6) edge[very thick] (bf)
		(aplBf1) edge node{$\hat{X}_{f,t,1}$} (NN2)
		(norm1) edge[draw=newCostCol,fill=newCostCol, very thick] (aplBf2)
		(norm2) edge[draw=newCostCol,fill=newCostCol, very thick] (aplBf3)
		(br7) edge[draw=newCostCol,fill=newCostCol, very thick] (aplBf2)
		(aplBf2) edge[dashed,draw=newCostCol,fill=newCostCol, very thick] (aplBf3)
		(aplBf2) edge[draw=newCostCol,fill=newCostCol] (\cSeven-0.4,\rTwo)
		(aplBf3) edge[draw=newCostCol,fill=newCostCol] (\cSeven-0.4,\rOne)
		(snr) edge[draw=newCostCol,fill=newCostCol] (lossNew)
		(\cTwo+0.49,\rThree+0.2) edge[draw=oldCostCol,fill=oldCostCol, very thick] (\cThree-0.5,\rThree+0.2)
		(\cTwo+0.49,\rThree-0.2) edge[draw=oldCostCol,fill=oldCostCol, very thick] (\cThree-0.5,\rThree-0.2)
		(br3) edge[draw=oldCostCol,fill=oldCostCol, very thick] (\cTwo-0.2,\rThree+0.3)
		(br4) edge[draw=oldCostCol,fill=oldCostCol, very thick] (\cTwo+0.2,\rThree+0.3)
		(br8) edge[draw=oldCostCol,fill=oldCostCol, very thick] (\cThree+0.5,\rThree+0.2)
		(br9) edge[draw=oldCostCol,fill=oldCostCol, very thick] (\cThree+0.5,\rThree-0.2)
		(bce) edge[draw=oldCostCol,fill=oldCostCol] (lossOld)
		(br10) edge[draw=werCostCol,fill=werCostCol] (ce)
		(p) edge[draw=werCostCol,fill=werCostCol] (ce)
		(ce) edge[draw=werCostCol,fill=werCostCol] (lossWer)
	;

	\draw[->,draw=newCostCol, very thick](br1) |- (norm2);
	\draw[->,draw=newCostCol, very thick](br2) |- (norm1);
	\draw[->, very thick](n) -| (sum);
	\draw[->, very thick](bf) -| node[above left]{$\w_f\hermite$} (aplBf1);
	\draw[->, very thick](\cFour-0.15,\rFour-0.3) |- (\cFife-0.4,\rThree-0.2);
	\draw[->, very thick](\cFour+0.15,\rFour-0.3) |- (\cFife-0.4,\rThree+0.2);

\end{tikzpicture} \medskip
	\caption{Block diagram of the training system of the described system with different training objectives in colors}
	\label{fig:blockschaltbild}
\end{figure}

Fig.~\ref{fig:blockschaltbild} gives an overview of the considered \gls{NN} supported acoustic beamformer for speech enhancement and the following \gls{asr} with all three used objective functions \cite{Heymann2015}.
The different colors mark the different objective functions which are used for training.
After both networks are trained only the black parts are needed to decode the spoken word sequence.
In \cite{Heymann2015,Heymann2016} the mask estimation \gls{NN} uses heuristic \glspl{ibm} as targets and \gls{bce} as the objective function. This is shown in blue.
Our new approach (shown in green) is to use the negative SNR of the beamformed signal as the objective function and propagate the gradient through the beamformer back to the mask estimation \gls{NN}.
The global approach is to use a common objective function for both front-end and back-end networks.
For the latter the cross entropy is used as the objective function which is well known for \gls{asr} tasks (here shown in purple) and the gradient is propagated back to the mask estimation \gls{NN} for optimization \wrt the cross entropy.
The reference target $p$ is the framewise label/state posterior probability or sequence probability (more information in \cite{heymann2017EndToEnd}).

The computation of the objective function $J$ involves complex-valued quantities although the objective is real-valued.
Thus complex-valued gradients are needed for error backpropagation, which are described in \cref{cha:complexbackprop}.

\section{Statistically optimal beamforming operations}
\label{sec:beamformers}
In this section different beamforming operations are presented. First the widely used \gls{mvdr} beamformer is described, followed by the \gls{gev} beamformer, which maximizes the \gls{snr} of the output signal.
\subsection{MVDR beamformer}
The \gls{mvdr} beamformer tries to minimizes the variance of the noise power without distorting the signal from the desired direction.
The coefficients can be calculated via
\begin{align}
	\w^{(\MVDR)}_{f} &= \frac{\phinn_f\inv \cdot \w^{(\PCA)}_{f} }{\w^{(\PCA)}_{f}{}\hermite \cdot \phinn_f\inv \cdot \w^{(\PCA)}_{f}},
\end{align}
where $\w^{(\PCA)}_{f} = \mathcal{P}(\phixx_f)$ is the eigenvector corresponding to the largest eigenvalue of $\phixx_f$ \cite{Haykin2002}. The gradients needed for forward and backward mode \gls{ad} of the \gls{mvdr} beamformer can be obtained with the examples found in \autoref{cha:blocks}.
\subsection{GEV beamformer}
The \gls{gev} beamforming vector \cite{warsitz2006gev} 
\begin{align}
	\w_f^{(\GEV)} = \argmax{\w_f}{\frac{\E\{|\w_f\hermite \vec{Y}_{f,t}|^2\mid \vec{N}_{f,t}=0\}}
						   {\E\{|\w_f\hermite \vec{Y}_{f,t}|^2\mid \vec{X}_{f,t}=0\}}}
				  = \argmax{\w_f}{\frac{\w_f\hermite\phixx_f\w_f}{\w_f\hermite\phinn_f\w_f}},
\end{align}
also known as max SNR beamfomrer, requires the solution of the generalized eigenvalue problem 
\begin{align}
	\phixx_f \vec{W}_f = \phinn_f \vec{W}_f \vecG{\Lambda}_f,
\end{align}
where $\w_f^{(\GEV)}$ is the eigenvector corresponding to the largest eigenvalue.
Under the assumption that the inverse of $\phinn_f$ exists, 
the generalized eigenvalue problem for hermitian matrices is equal to a normal eigenvalue problem 
\begin{align}
	\phinn_f\inv \phixx_f \vec{W}_f = \vec{W}_f \vecG{\Lambda}_f, \\
	\vecG{\Phi}_f \vec{W}_f = \vec{W}_f \vecG{\Lambda}_f,
\end{align}
where usually $\vecG{\Phi}_f$ is a non hermitian matrix.
This formulation is required to obtain the gradient, which is explained in the following chapters.


\chapter{Complex-valued backpropagation}
\label{cha:complexbackprop}
In the forward step of a \gls{NN} an objective from the input data and the network parameters, which should be optimized, is calculated.
Then the backpropagation algorithm uses the chain rule to calculate the gradients of the objective function \wrt all adjustable parameters.
With these gradients an optimizer (e.g. gradient decent) updates the parameters.
The real-valued gradient calculation rules, especially the chain rule
\begin{align}
	\der{J}{x} = \der{J}{\omega}\der{\omega}{x},
\end{align}
are well known. 
This gets more difficult if complex numbers are involved in the forward step.
There is a definition for complex differentiation: 
A function is complex differentiable (holomorph) if the real and imaginary part satisfy the Cauchy–Riemann equations \cite[Chapter 13]{kreyszig2010math}. 
Since many relevant complex function are not holomorph (e.g. absolute value, real/imaginary part and conjugation) we dismiss this approach. 

An alternative definition for the complex derivative is described in \cite{Brandwood1983} the {\em Wirtinger calculus}
\begin{align}
	\der{J}{z\conj} &= \frac{1}{2} \left(\der{J}{x} + \j \der{J}{y}\right),\label{eq:Wirtinger} \\
	\der{J}{z} &= \frac{1}{2} \left(\der{J}{x} - \j \der{J}{y}\right),\label{eq:Wirtinger2}
\end{align}
where $J$ is the \rv objective, $x$ and $y$ are real numbers and $\der{J}{x}$ is the partial derivative of $J$ \wrt $x$.
This definition is well known for complex optimization problems. 
It can be interpreted as a pure real-valued problem, because the objective function is real-valued and all complex operations can be separated in real operations. 
Nevertheless in this report we use the {\em Wirtinger calculus} to derive rules for calculating complex gradients instead of defining the forward step via real-valued operations.

Since an objective function has to be real valued it is easy to switch between \eqref{eq:Wirtinger} and \eqref{eq:Wirtinger2}
\begin{align}
\der{J}{z\conj} &= \left(\der{J}{z}\right)\conj.
\end{align}
Although both definitions can be used we mainly focus on \eqref{eq:Wirtinger} here.

\section{Chain rule}
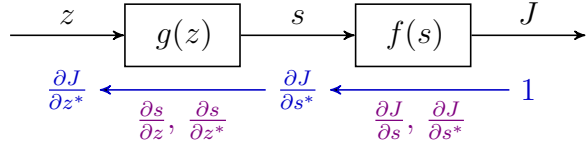
\begin{figure}[H]
	\centering
	\definecolor{backpropcol}{RGB}{0,0,200}
	\definecolor{derivativecol}{RGB}{128,0,128}

\begin{tikzpicture}
	 [>=stealth',semithick,auto,node distance = 3cm,
	  	block/.style={rectangle,draw,text centered, text width=3em, minimum height=2em},
	   	branch/.style = {circle,inner sep=0pt,minimum size=1mm,fill=black,draw=black},
	   	clear/.style = {minimum height=2em, text centered,text=backpropcol}
	 ]

	 \node[text width = 0em, text centered](z){};
	 \node[block, right = 1.5cm of z](g){$g(z)$};
	 \node[block, right = 1.5cm of g](f){$f(s)$};
	 \node[text centered, right = 1.5cm of f](J){};

	 \path [{}-{>}]	
	 	(z) edge node(z){$z$} (g)
	 	(g) edge node(s){$s$} (f)
	 	(f) edge node(j){$J$}(J);

	 \node[clear, below = 0.3cm of s](dJds){$\der{J}{s\conj}$};
	 \node[clear, below = 0.3cm of z](dJdz){$\der{J}{z\conj}$};
	 \node[clear, below = 0.3cm of j](one){1};
	 
	 \path [{}-{>},draw=backpropcol]	
	  	(one) edge node[below,text=derivativecol]{$\der{J}{s}$, $\der{J}{s\conj}$} (dJds)
	  	(dJds) edge node[below,text=derivativecol]{$\der{s}{z}$, $\der{s}{z\conj}$} (dJdz);
\end{tikzpicture}%
 \medskip
	\caption[Visualization of the chain rule: forward step, backward gradient and intermediate derivatives]{Visualization of the chain rule: forward step, {\color{backpropcol} backward gradient} and {\color{derivativecol} intermediate derivatives}}
	\label{fig:chain_rule}
\end{figure}
The most important rule in this work is the chain rule. 
This rule describes the relation between gradients,
when a function is applied to the output of another. 
Therefore, we start with a chain of two functions (see \cref{fig:chain_rule})
\begin{align}
\begin{aligned}
J = f(s) &\quad \text{where }  s = \sigma + \j \omega,\\
s = g(z) &\quad \text{where } z = x + \j y, \\
s, z \in \Complex, &\quad J, \sigma, \omega, x, y \in \Real.
\end{aligned}
\end{align}
Combining the real-valued chain rule \cite[p. 393]{kreyszig2010math}
\begin{align}
\der{J}{x} = \der{J}{\sigma} \der{\sigma}{x} + \der{J}{\omega} \der{\omega}{x}
\end{align}
and \eqref{eq:Wirtinger} leads to 
\begin{align}
\der{J}{z\conj} &=\frac{1}{2} \left(\der{J}{\sigma} \der{\sigma}{x} + \der{J}{\omega} \der{\omega}{x}+ \j \der{J}{\sigma} \der{\sigma}{y} + \j \der{J}{\omega} \der{\omega}{y}\right)\nonumber\\
&= \dots\text{\footnotemark} = \left(\der{J}{s\conj}\right)\conj \der{s}{z\conj} + \der{J}{s\conj} \left(\der{s}{z}\right)\conj, \label{eq:chainRule}
\end{align}
\footnotetext{A detailed calculation can be found in \cref{appendix:proof:complex_valued_chain_rule}}
where
\begin{align}
\begin{aligned}
\der{s}{z\conj} &= \frac{1}{2} \left(\der{s}{x} + \j \der{s}{y}\right) = \frac{1}{2} \left(\der{\sigma}{x} + j \der{\omega}{x} + \j \der{\sigma}{y} - \der{\omega}{y}\right), \\
\der{s}{z} &= \frac{1}{2} \left(\der{s}{x} - \j \der{s}{y}\right) = \frac{1}{2} \left(\der{\sigma}{x} + j \der{\omega}{x} - \j \der{\sigma}{y} + \der{\omega}{y}\right). \label{eq:wirtinger_for_s}
\end{aligned}
\end{align}
This solution for the chain rule is the same as in \cite{kreutz2009WirtingerChainRule,amin2011wirtinger,bouboulis2011WirtingerChainRule,wiki:wirtinger}

\section{Reuse of a variable}
\label{sec:reuse_variable}

Sometimes a variable is used more than one time for the calculation of the objective function. 
Thus this variable receives more than one gradient which are accumulated like in the real-valued case.

With this rule we can extend the scalar chain rule to tensors with arbitrary dimensions. 
When $\vec{S}$ and $\vec{Z}$ are tensors
\begin{align}
\vec{S} &\in \Complex^{A_1 \times \dots \times A_N} ~,~~
\vec{Z} \in \Complex^{B_1 \times \dots \times B_M} \nonumber
\end{align}	
the chain rule can be denoted elementwise 
\begin{align}
\der{J}{Z_{b_1, \dots, b_M}\conj} &= \sum_{a_1, \dots, a_N} \left(\der{J}{S_{a_1, \dots, a_N}\conj}\right)\conj \der{S_{a_1, \dots, a_N}}{Z_{b_1, \dots, b_M}\conj} + \der{J}{S_{a_1, \dots, a_N}\conj} \left(\der{S_{a_1, \dots, a_N}}{Z_{b_1, \dots, b_M}}\right)\conj \label{eq:chain_rule_tensor},
\end{align}	
where $a_n$ and $b_m$ ($n \in \{1, \dots, N\}$ and $m \in \{1, \dots, M\}$) indicate the elements in the $n$-th and $m$-th dimension, respectively.

A special case of this chain rule, where only $\vec{S}$ is a two dimensional matrix and $z$ is scalar, can be written as
\begin{align}
\der{J}{z\conj} &= \sum_{n, m}\left(\left(\der{J}{S_{n, m}\conj}\right)\conj \der{S_{n, m}}{z\conj} + \der{J}{S_{n, m}\conj} \left(\der{S_{n, m}}{z}\right)\conj\right) \nonumber \\
&= \trace{\left(\der{J}{\vec{S}\conj}\right)\hermite \der{\vec{S}}{z\conj} + \left(\der{J}{\vec{S}\conj}\right)\transposed \left(\der{\vec{S}}{z}\right)\conj}, \label{eq:complChainRule}
\end{align}
where $\trace{\cdot}$ is the sum of the diagonal elements, $\der{J}{\vec{S}\conj}$ the derivative of $J$ \wrt each element of $\vec{S}$ and $\der{\vec{S}}{z\conj}$ the derivative of each element of $\vec{S}$.
Both gradient matrices have the same shape as $\vec{S}$.
This notation using the $\trace{\cdot}$ operation is only useful for symbolic calculations in an algorithm explained later.

To compare this chain rule with \cite[Eq. 7]{Drude2016CVNN}, $\vec{s}$ and $\vec{z}$ have to be vectors
\begin{align}
\der{J}{\vec{z}\conj} 
&= \left(\left(\der{J}{\vec{s}\conj}\right)\hermite \der{\vec{s}}{\vec{z}\conj} 
+ \left(\der{J}{\vec{s}\conj}\right)\transposed \left(\der{\vec{s}}{\vec{z}}\right)\conj\right)\transposed \label{eq:chain_rule_vector},\\
\der{\vec{s}}{\vec{z}\conj} 
&= \left[\begin{array}{ccc}
\der{s_1}{z_1\conj} & \der{s_1}{z_2\conj} & \dots\\
\der{g_2}{z_1\conj} & \der{s_2}{z_2\conj} & \\
\vdots &  & \ddots
\end{array}\right]. \label{eq:def_der_vec_wrt_vec}
\end{align}

\section{Real-valued input}
\label{sec:realValuedInput}

One important case of a complex operation is a real-valued input. 
This means, that the imaginary part $\Im\{z\} = y$ is zero, so $s$ is independent of it. 
Therefore the gradient \wrt $y$ has to be zero and \eqref{eq:chainRule} simplifies to
\begin{align}
\der{J}{z\conj} 
&= \frac{1}{2} \left(\der{J}{\sigma} \der{\sigma}{x} + \der{J}{\omega} \der{\omega}{x}+ \j \der{J}{\sigma} 0 + \j \der{J}{\omega} 0\right) \nonumber\\
&= \frac{1}{2} \left(\der{J}{\sigma} \der{\sigma}{x} + \der{J}{\omega} \der{\omega}{x}\right) \nonumber\\
&= \Re\left\{\left(\der{J}{s\conj}\right)\conj \der{s}{z\conj} + \der{J}{s\conj} \left(\der{s}{z}\right)\conj\right\},
\end{align}
where $\Re\left\{\cdot\right\}$ denotes the real part.
Note that the real part operator has to be applied to the gradient if the input is real-valued,
because in general the gradient is \cv.

\section{Gradient descent}

Optimizing the parameters of a \gls{NN} can be done with gradient descent.
For the real-valued case this is done via
\begin{align}
	x^{(\mathrm{new})} &= x^{(\mathrm{old})} - \frac{\mu}{2} \der{J}{x} \bigg|_{x = x^{(\mathrm{old})}}, \\
	y^{(\mathrm{new})} &= y^{(\mathrm{old})} - \frac{\mu}{2} \der{J}{y} \bigg|_{y = y^{(\mathrm{old})}},
\end{align}
where $\mu / 2$ is the stepsize.
Combining both equations with the definition of complex numbers and the {\em Wirtinger calculus} \eqref{eq:Wirtinger}
\begin{align}
	z^{(\mathrm{new})} &= x^{(\mathrm{new})} + \j y^{(\mathrm{new})} \nonumber\\
	&= x^{(\mathrm{old})} - \frac{\mu}{2} \der{J}{x} \bigg|_{x = x^{(\mathrm{old})}} + \j y^{(\mathrm{old})} - \frac{\mu}{2} \der{J}{y} \bigg|_{y = y^{(\mathrm{old})}} \nonumber\\
	&= x^{(\mathrm{old})} + \j y^{(\mathrm{old})} - \frac{\mu}{2} \left(\der{J}{x} \bigg|_{x = x^{(\mathrm{old})}} + \j \der{J}{y} \bigg|_{y = y^{(\mathrm{old})}}\right) \nonumber\\
	&= z^{(\mathrm{old})} - \mu \der{J}{z\conj} \bigg|_{z = z^{(\mathrm{old})}} \label{eq:gradDesc1} \\
	&= z^{(\mathrm{old})} - \mu \left(\der{J}{z}\right)\conj \bigg|_{z = z^{(\mathrm{old})}} \label{eq:gradDesc2}
\end{align}
results in a solution for the complex gradient descent.
Note, that there are two possibilities to calculate the gradient.
The first one, shown in \eqref{eq:gradDesc2}, is to determine $\der{J}{z}$ and use its conjugate complex as the update like it is done in the framework \cite{Maclaurin2015autograd}.
Otherwise $\der{J}{z\conj}$ is determined and the same update rule as in the real-valued case is used, denoted in \eqref{eq:gradDesc1}.

\section{Derivation of the backward gradient}
The gradient of a function is needed for the training of a \gls{NN} with the backpropagation method.
Here this is done via \gls{ad} where a function is segmented into a chain of elementary operations.
The input gradient of each elementary operation is calculated with the chain rule from the output gradient.
This is done to get the gradient \wrt the target quantities.
This method is called building block concept. 

The idea is to build a computational graph for the complete network.
Beginning at the end, each block receives the output gradient 
and calculates the input gradient using the chain rule until the gradient of each weight is known. 
The building block system can also be used to analytically combine intermediate blocks to one block with a numerically more stable gradient and an improved execution time.


In this section different strategies to derive the gradient of a function are presented.
First\inred{,} two analytical derivations are described followed by a numerical approach.
The numerical approach is used to verify the implementation of the analytical gradients.
These strategies are only two suggestions to get a solution.
There may be others which are better suited for different tasks.

\subsection{Conventional derivation} 
\label{sec:der:drude}

Following the ideas from \cite{Drude2016CVNN} to derive the gradient of a complex valued function is intuitive.
After reformulating the function to a scalar version it is possible to calculate the derivative of each output \wrt each input ($\der{s_n}{z_m}$,  $\der{s_n}{z_m\conj}$).
Sometimes it is helpful to use the real-valued derivatives ($\der{\sigma}{x},\der{\sigma}{y},\der{\omega}{x},\der{\omega}{y}$).
Now the derivatives can be obtained by inserting the derivatives in \eqref{eq:chain_rule_tensor} or \eqref{eq:chain_rule_vector} in the case of vectors.

This approach, summarized in \autoref{alg:Drude}, to calculate the derivative is useful when the function of interest describes scalar operation.
However a disadvantage is that a non implicit function is required, i.e. the target has to be alone on one side of the equal sign.
The algorithm in the following section eliminates this drawback.

\begin{algorithm}[H]
	\caption{Conventional derivation}
	\label{alg:Drude}
	\begin{algorithmic}[1]
		\State Calculate the derivative of each output \wrt each input: $\der{s_n}{z_m}$, $\der{s_n}{z_m\conj}$
		\State Insert in \eqref{eq:chain_rule_tensor} or \eqref{eq:chain_rule_vector} in the case of vectors
	\end{algorithmic}
\end{algorithm}

\subsection{Derivation over forward derivative} 
\label{sec:der:giles}

This complex-valued approach as a generalization of \cite{Giles2008} differs from the previous section.
Here the initial equation is a matrix operation $\vec{C} = \vec{h}(\vec{A}, \vec{B})$, where \vec{C} is the output and \vec{A} and \vec{B} are the inputs.
The equation is reformulated to a simple characteristic equation which describes the dependencies between the matrices, e.g. $\vec{A}\vec{C} = \vec{B}$ is a simpler formulation as $\vec{C} = \vec{A}\inv\vec{B}$ (\cref{subsec:matrixProdInv}).
Now both sides of the equation are derived \wrt $\zeta \in \{z, z\conj\}$ where $\zeta$ is assumed as the input variable before this operation.
With the reformulation of this equation to $\derz{\vec{C}} = \dots$ the forward mode \gls{ad} is obtained \cite{Giles2008} which is not required for the backpropagation algorithm.\footnote{
	The forward mode \gls{ad} calculates the gradient by starting at the input and propagate through all calculations.
	Therefore it requires for each parameter an independent calculation of the gradient. 
	In many applications this is too expensive and the forward mode \gls{ad} $(\der{\vec{C}}{z\conj})$ is used as an analytical helper to receive the backward mode \gls{ad} $(\der{J}{\vec{A}\conj})$. 
	Note: In forward mode the derivation of the output $\vec{C}$ \wrt the input $z\conj$ and in backward mode the derivation of the objective function \wrt the input \vec{A} is required.}
To obtain the backward mode \gls{ad}, needed for optimization, the forward mode AD is inserted in \eqref{eq:complChainRule}
\begin{align}
\der{J}{z\conj} &= \trace{\left(\der{J}{\vec{C}\conj}\right)\hermite \der{\vec{C}}{z\conj} + \left(\der{J}{\vec{C}\conj}\right)\transposed \left(\der{\vec{C}}{z}\right)\conj} \label{eq:giles:trace_c}
\end{align}
and rearranged to 
\begin{align}
\der{J}{z\conj} &= \trace{\left(\dots\right)\hermite \der{\vec{A}}{z\conj} + \left(\dots\right)\transposed \left(\der{\vec{A}}{z}\right)\conj + \left(\dots\right)\hermite \der{\vec{B}}{z\conj} + \left(\dots\right)\transposed \left(\der{\vec{B}}{z}\right)\conj}.
\end{align}
A comparison of the coefficients with
\begin{align}
\der{J}{z\conj} &= \trace{\left(\grad{\vec{A}}\right)\hermite \der{\vec{A}}{z\conj} + \left(\grad{\vec{A}}\right)\transposed \left(\der{\vec{A}}{z}\right)\conj + \left(\grad{\vec{B}}\right)\hermite \der{\vec{B}}{z\conj} + \left(\grad{\vec{B}}\right)\transposed \left(\der{\vec{B}}{z}\right)\conj} \label{eq:giles:trace_a_b}
\end{align}
results in the backward mode AD $\grad{\vec{A}}$ and $\grad{\vec{B}}$. 
Therefore the use of the trace identities like 
\begin{align}
\trace{A\transposed}&=\trace{A}, \\
\trace{AB}&=\trace{BA}, \\
\trace{A+B}&=\trace{A}+\trace{B}
\end{align}
is beneficial.

With this concept, described in \autoref{alg:Giles}, it is possible to find more solutions than with the previous one.
For example the backward gradient of the inverse $\vec{C} = \vec{A}\inv$ (\cref{subsec:matrixInv}) can be obtained which can not be decomposed into a trivial element-wise operation.

\begin{algorithm}[H]
	\caption{Derivation over forward derivative}
	\label{alg:Giles}
	\newcommand{\IndState}{\Statex\hspace{\algorithmicindent}}  
	\begin{algorithmic}[1]
		\State $\vec{C} = \vec{h}(\vec{A}, \vec{B})$ with matrices $\vec{A}$, $\vec{B}$ and $\vec{C}$ 
		\State Find a simple representation of the equation (Implicit functions are allowed)
		\State Compute the derivative of both sides w.r.t. $z\conj$ and $z$
		\State Reformulate and insert $\der{\vec{C}}{z\conj}$ and $\der{\vec{C}}{z}$ in \eqref{eq:complChainRule} 
		\State Use trace identities to reformulate
		\State Compare the coefficients with: 
		\IndState $\der{J}{z*} = \trace{\left(\der{J}{\vec{A}\conj}\right)\hermite \der{\vec{A}}{z\conj} + \left(\der{J}{\vec{A}\conj}\right)\transposed \left(\der{\vec{A}}{z}\right)\conj} + \trace{\left(\der{J}{\vec{B}\conj}\right)\hermite \der{\vec{B}}{z\conj} + \left(\der{J}{\vec{B}\conj}\right)\transposed \left(\der{\vec{B}}{z}\right)\conj}$
		\State The comparison of $\left(\der{J}{\vec{A}\conj}\right)\hermite$ and $\left(\der{J}{\vec{A}\conj}\right)\transposed$ should result in the same solution for $\der{J}{\vec{A}\conj}$
	\end{algorithmic}
\end{algorithm}	

\subsection[Numeric derivative]{Numeric derivative \cite{chainer_learningsys2015}}
\label{sec:numeric:derivative}

This concept can be used to verify the analytical solutions.
The definition of the derivative is used to calculate an approximation for the derivative.
Therefore only the forward step needs to be known so that the algorithm can approximate the derivative for a specific input value.
The \cref{alg:numerical} is very slow but a good way to verify the correctness of the analytical solution.
The operator $\mathrm{vec}$ in the \cref{alg:numerical} is the vectorize or flattern operator. 
It takes the input with arbitrary numbers of dimensions and concatenate all elements to a vector. 
The new order of the elements is not necessary in this example.
\begin{algorithm}[H]
	\caption{Numerical derivative}
	\label{alg:numerical}
	\newcommand{\lastState}{$\der{J}{z\conj} 
		\approx \left(\der{J}{\vec{s}\conj}\right)\hermite 
		\frac{1}{2}\left(\der{\vec{s}}{x} + \j \der{\vec{s}}{y}\right) 
		+ \left(\der{J}{\vec{s}\conj}\right)\transposed 
		\left(\frac{1}{2}\left(\der{\vec{s}}{x} - \j \der{\vec{s}}{y}\right)\right)\conj$}
	\begin{algorithmic}[1]
		\State $\vec{s} = \vec{g}(z) \gets \mathrm{vec}\{\vec{g}(z)\}$ \vphantom{$\left(\frac{\der{\vec{s}}{x} - \j \der{\vec{s}}{y}}{2}\right)\conj$}
		\State $\der{\vec{s}}{x} \approx \frac{\vec{g}(x+\varepsilon+\j y) - \vec{g}(x-\varepsilon+\j y)}{2\varepsilon}$\vphantom{\lastState}
		\State $\der{\vec{s}}{y} \approx \frac{\vec{g}(x+\j \left(y+\varepsilon\right)) - \vec{g}(x+\j \left(y-\varepsilon\right))}{2\varepsilon}$\vphantom{\lastState}
		\State \lastState
	\end{algorithmic}
\end{algorithm}


\chapter{Blocks}
\label{cha:blocks}

In this chapter we present some examples for the \cv backpropagation.
Starting with some fundamental elementwise operations it finally ends with the derivation of the eigenvalue problem.
All described operations are tested against a numerical gradient.
Many functions were developed for \cite{Drude2016CVNN}, but they are not described there.
The derivation of the matrix operations are inspired by their \rv counterparts in \cite{Giles2008}.
Also the solutions are compared with their \rv implementation in \cite{chainer_learningsys2015} and \cv implementation in \cite{Maclaurin2015autograd} in case they exists.

\section{Identity}

\newcommand{\gradIdentity}{\grad{s}}

The simplest possible example is the identity 
\begin{align}
	s = g(z) = z.
\end{align}
With the derivatives
\begin{align}
	\der{s}{z} &= 1 \qquad
	\text{and}\qquad	\der{s}{z\conj} = 0 \nonumber
\end{align}
and \eqref{eq:chainRule} the gradient is
\begin{align}
	\der{J}{z\conj} &= \left(\der{J}{s\conj}\right)\conj 0 + \der{J}{s\conj} (1)\conj 
	= \gradIdentity. \label{eq:gradIdentity}
\end{align}
As expected the identity does not change the gradient.

\section{Conjugation}
\label{sec:block:conj}
\newcommand{\gradConj}{\left(\grad{s}\right)\conj}

An important operation for complex numbers is the conjugation 
\begin{align}
	s = g(z) = z\conj.
\end{align}
The derivatives are again as trivial as for the identity 
\begin{align}
\der{s}{z} &= 0 \qquad
\text{and}\qquad	\der{s}{z\conj} = 1 \nonumber
\end{align}
and with \eqref{eq:chainRule} the gradient is
\begin{align}
\der{J}{z\conj} &= \left(\der{J}{s\conj}\right)\conj 1 + \der{J}{s\conj} (0)\conj 
= \gradConj. \label{eq:gradConj}
\end{align}
In case that $z$ is \rv the gradient $\der{J}{s\conj}$ is \rv, too, and the conjugation merges into the identity as expected.

\section{Addition}
\label{sec:block:add}
\newcommand{\gradAdditionOne}{\grad{s}}
\newcommand{\gradAdditionTwo}{\grad{s}}

Let 
\begin{align}
	s = g(z_1, z_2) = z_1 + z_2,
\end{align}
then the gradients are 
\begin{align}
\grad{z_1} &= \gradAdditionOne \qquad \text{and} \qquad \grad{z_2} = \gradAdditionTwo \label{eq:gradAdd}
\end{align}
like the \rv gradient.

\section{Multiplication}
\label{sec:block:mul}
\newcommand{\gradMultiplicationOne}{\grad{s} z_2\conj}
\newcommand{\gradMultiplicationTwo}{\grad{s} z_1\conj}

In the case  that 
\begin{align}
	s = g(z_1, z_2) = z_1 \cdot z_2
\end{align}
the derivatives are
\begin{align}
\der{s}{z_1} &= z_2 \qquad
\text{and}\qquad	\der{s}{z_1\conj} = 0, \nonumber\\
\der{s}{z_2} &= z_1 \qquad
\text{and}\qquad	\der{s}{z_2\conj} = 0 \nonumber
\end{align}
and with \eqref{eq:chainRule} the gradients are
\begin{align}
\der{J}{z_1\conj} 
&= \left(\der{J}{s\conj}\right)\conj 0 + \der{J}{s\conj} z_2\conj
= \gradMultiplicationOne, \\
\grad{z_2} 
&= \left(\grad{s}\right)\conj 0 + \grad{s} z_1\conj \label{eq:gradMul}
= \gradMultiplicationTwo.
\end{align}

\section{Exponentiation}
\label{sec:block:pow}
\newcommand{\gradPow}{\grad{s} n\left(z\conj\right)^{n-1}}

The square 
\begin{align}
	s = g(z) = z^2
\end{align}
is a special case of the Multiplication ($z=z_1=z_2$). Therefore the gradient is
\begin{align}
\der{J}{z\conj} 
&= \der{J}{s\conj} z_1\conj + \der{J}{s\conj} z_2\conj = 2\der{J}{s\conj} z\conj. \nonumber
\end{align}
For the generalization $s = g(z) = z^n$ we split the equation in $s = z_1 z_2$ $n$-times, where $z_1 = z$ and $z_2=z^{n-1}$. With the gradient of the multiplication
\begin{align}
	\grad{z_1} = \der{J}{s\conj} z_2\conj = \der{J}{s\conj} \left(z\conj\right)^{n-1} \nonumber
\end{align}
and the knowledge of \cref{sec:reuse_variable} to accumulate all $n$ gradients, we get
\begin{align}
\grad{z} 
&= \gradPow. \label{eq:gradPow}
\end{align}

\section{Division}
\label{sec:block:div}
\newcommand{\gradDivisionOne}{\grad{s}\frac{1}{z_2\conj}}
\newcommand{\gradDivisionTwo}{\grad{s}\frac{-z_1\conj}{\left(z_2\conj\right)^2}}

The derivation of the division is more complicated. Therefore define 
\begin{align}
	s = g(z_1, z_2) = z_1 / z_2 = z_1 z_3
\end{align}
where $z_3 = 1 / z_2$. We know from multiplication
\begin{align}
\grad{z_1} = \grad{s}z_3\conj = \gradDivisionOne
\end{align}
and
\begin{align}
\grad{z_3}=\grad{s}z_1\conj . \nonumber
\end{align}
From exponentiation with $n=-1$ we know that
\begin{align}
\grad{z_2} = \grad{z_3}\left(-1\right)\left(z_2\conj\right)^{-2}\nonumber
\end{align}
Therefore the gradient is
\begin{align}
\grad{z_2} = \gradDivisionTwo. \label{eq:gradDiv}
\end{align}

\section{Absolute value}
\label{sec:block:abs}
\newcommand{\gradAbsolute}{\grad{s} \e^{\j \varphi}}

To derive the gradient of the absolute value 
\begin{align}
	s = g(z) = \left|z\right| \in \mathbb{R}
\end{align}
we use the real derivatives with $z = x + \j y = \left|z\right| \e^{\j\varphi}$
\begin{align}
	\der{s}{x} &= \der{\sqrt{x^2 + y^2}}{x} = \frac{x}{s} \qquad \text{and} \qquad 
	\der{s}{y} = \der{\sqrt{x^2 + y^2}}{y} = \frac{y}{s} \nonumber
\end{align}
and insert them in the definition \eqref{eq:wirtinger_for_s}
\begin{align}
\der{s}{z\conj} &= \frac{1}{2}\left(\der{s}{x} + \j \der{s}{y}\right) = \frac{1}{2}\frac{z}{s} = \frac{1}{2} \e^{\j \varphi}, \nonumber\\
\der{s}{z} &= \frac{1}{2}\left(\der{s}{x} - \j \der{s}{y}\right) = \frac{1}{2}\frac{z\conj}{s} = \frac{1}{2} \e^{-\j \varphi}. \nonumber
\end{align}
Using \eqref{eq:chainRule} leads to 
\begin{align}
\der{J}{z\conj} 
&= \underbrace{\left(\der{J}{s\conj}\right)\conj}_{(\cdot)\in \mathbb{R}} \frac{1}{2} \e^{\j \varphi} + \der{J}{s\conj} \left(\frac{1}{2} \e^{-\j \varphi}\right)\conj
= \gradAbsolute, \label{eq:gradAbs}
\end{align}
where $\der{J}{s\conj}$ is \rv, because it is the gradient of a \rv number.

\section{Phase factor} 
\newcommand{\gradExpAngle}{\grad{s} \frac{1}{\left|z\right|} - \Re\left\{\grad{s} \frac{1}{z}\right\} \e^{\j \varphi}}

For the operation 
\begin{align}
	s = g(z) = \e^{\j\varphi} = \frac{z}{\left|z\right|}
\end{align}
(sometimes called the generalization of the sign function to complex numbers) we want to demonstrate two derivations. One way over the real derivatives (left column) and one by using the chain rule with known derivatives (right column). 


\begin{paracol}{2}
	\setlength{\columnseprule}{.4pt}
	\setlength{\parskip}{0pt} 
	\begin{leftcolumn}
		For the \rv derivative we insert $z = x + \j y$ in the operation
		\begin{align}
			s = \e^{\j\varphi} = \frac{x + \j y}{\sqrt{x^2 + y^2}} \nonumber
		\end{align}
		and calculate the derivative \wrt $x$
		\begin{align}
			\der{s}{x} &= \frac{1\sqrt{x^2+y^2}-(x+\j y)\frac{x}{\sqrt{x^2+y^2}}}{x^2+y^2} \nonumber\\
			&= \frac{\left|z\right| - x \frac{z}{\left|z\right|}}{\left|z\right|^2}.\nonumber
		\end{align}
		The calculation \wrt $y$ is analogue
		\begin{align}
			\der{s}{y} &= \frac{j\left|z\right| - y \frac{z}{\left|z\right|}}{\left|z\right|^2}. \nonumber
		\end{align}
		Insert both in the definition \eqref{eq:wirtinger_for_s}
		\begin{align}
			\der{s}{z\conj} &= \frac{1}{2}\left(\frac{\left|z\right|-x\frac{z}{\left|z\right|}+\j\cdot\j \left|z\right| - \j y \frac{z}{\left|z\right|}}{\left|z\right|^2}\right) \nonumber\\
			&= \frac{-z}{2\left|z\right|^3} \left(x +\j y\right) \nonumber\\
			&= \frac{-1}{2\left|z\right|}\e^{\j2\varphi}, \nonumber\\
			\der{s}{z} &= \frac{1}{2}\left(\frac{\left|z\right|-x\frac{z}{\left|z\right|}-\j\cdot\j \left|z\right| + \j y \frac{z}{\left|z\right|}}{\left|z\right|^2}\right) \nonumber\\
			&= \frac{1}{\left|z\right|}-\frac{z}{2\left|z\right|^3} \left(x -\j y\right) \nonumber\\
			&= \frac{1}{\left|z\right|}-\frac{1}{2\left|z\right|}
			= \frac{1}{2\left|z\right|}\nonumber
		\end{align}
		and \eqref{eq:chainRule} results in the solution: 
		\begin{align}
			\grad{z} &= \left(\grad{s}\right)\conj \frac{-1}{2\left|z\right|}\e^{\j2\varphi} + \grad{s} \left(\frac{1}{2\left|z\right|}\right)\conj . \nonumber
		\end{align}
	
	\end{leftcolumn}
	\begin{rightcolumn}
	
		Using the chain rules requires to introduce intermediate variables
		\begin{align}
			s = \e^{\j\varphi} = \frac{z}{\left|z\right|} =: \frac{z_1}{\left|z_3\right|} =: \frac{z_1}{z_2},\\
			z_2 \in \mathbb{R}.
		\end{align}
		From the division we know
		\begin{align}
			\grad{z_1} = \grad{s} \frac{1}{z_2\conj}, \nonumber \\
			\grad{z_2} = \Re\left\{\grad{s} \frac{-z_1\conj}{(z_2\conj)^2}\right\}, \nonumber
		\end{align}
		because $z_2$ is \rv $\grad{z_2}$ has to be forced to be \rv (\cref{sec:realValuedInput}), too.
		Next use the absolute value
		\begin{align}
			\grad{z_3} &= \grad{z_2} \e^{\j \varphi} \nonumber\\
			&= \Re\left\{\grad{s} \frac{-z_1\conj}{(z_2\conj)^2}\right\} \e^{\j \varphi}\nonumber
		\end{align}
		and combine the gradients
		\begin{align}
			\grad{z} &= \grad{z_1} + \grad{z_3} \nonumber\\
			&= \grad{s} \frac{1}{z_2\conj} + \Re\left\{\grad{s} \frac{-z_1\conj}{(z_2\conj)^2}\right\} \e^{\j \varphi} \nonumber \\
			&= \grad{s} \frac{1}{\left|z\right|} + \Re\left\{\grad{s} \frac{-z\conj}{\left|z\right|^2}\right\} \e^{\j \varphi} \nonumber\\
			&= \gradExpAngle. \label{eq:gradExpAngle}
		\end{align}
		For comparison with the \rv solution reformulate the equation
		\begin{align}
			\grad{z} &= \grad{s} \frac{1}{\left|z\right|} - \frac{1}{2} \left(\grad{s} \frac{1}{z} + \left(\grad{s} \frac{1}{z}\right)\conj\right) \e^{\j \varphi}	\nonumber\\
			&= \frac{1}{2} \grad{s} \frac{1}{\left|z\right|} - \frac{1}{2} \left(\grad{s}\right)\conj \frac{\e^{\j 2\varphi}}{\left|z\right|} . \nonumber
		\end{align}	
	\end{rightcolumn}
\end{paracol}

As expected both derivations results in the same solution. 
For the \rv way it is easy to derive the gradient, because all derivation rules are well known and the solution must only be inserted in the definitions. 
The complex chain rules has to be used more carefully, because when an intermediate variable is \rv the gradient must be forced to be \rv, too.
\pagebreak

\section{Real/Imaginary part}
\newcommand{\gradReal}{\grad{s}}
\newcommand{\gradImag}{\j\grad{s}}

The calculations for the real and imaginary part are similar.
\begin{paracol}{2}
	\setlength{\parskip}{0pt} 
	\setlength{\columnseprule}{.4pt}
	\begin{leftcolumn}
		The starting equation for the real part is
		\begin{align}
			s = g(z) = \Re\{z\} = \Re\{x+\j y\} = x
		\end{align}
		The derivatives \wrt $x$ and $y$ are
		\begin{align}
			\der{s}{x} &= \der{x}{x} = 1 \quad \text{and} \quad 
			\der{s}{y} = \der{x}{y} = 0.  \nonumber
		\end{align}
		Inserting in \eqref{eq:wirtinger_for_s}
		\begin{align}
			\der{s}{z\conj} &= \frac{1}{2}\left(\der{s}{x} + \j \der{s}{y}\right) = \frac{1}{2},  \nonumber\\
			\der{s}{z} &= \frac{1}{2}\left(\der{s}{x} - \j \der{s}{y}\right) = \frac{1}{2},  \nonumber
		\end{align}
		and \eqref{eq:chainRule}
		\begin{align}
			\der{J}{z\conj} 
			&= \underbrace{\left(\der{J}{s\conj}\right)\conj}_{(\cdot)\in \mathbb{R}} \frac{1}{2} + \der{J}{s\conj} \left(\frac{1}{2} \right)\conj \nonumber \\
			&= \gradReal \label{eq:gradReal}
		\end{align}
		yields the gradient.
	\end{leftcolumn}
	\begin{rightcolumn}
		
		The starting equation for the imaginary part is
		\begin{align}
			s = g(z) = \Im\{z\} = \Im\{x+\j y\} = y
		\end{align}
		The derivatives \wrt $x$ and $y$ are
		\begin{align}
			\der{s}{x} &= \der{y}{x} = 0 \quad \text{and} \quad 
			\der{s}{y} = \der{y}{y} = 1.  \nonumber
		\end{align}
		Inserting in \eqref{eq:wirtinger_for_s}
		\begin{align}
			\der{s}{z\conj} &= \frac{1}{2}\left(\der{s}{x} + \j \der{s}{y}\right) = \j \frac{1}{2},  \nonumber\\
			\der{s}{z} &= \frac{1}{2}\left(\der{s}{x} - \j \der{s}{y}\right) =  -\j \frac{1}{2},  \nonumber
		\end{align}	
		and \eqref{eq:chainRule}
		\begin{align}
			\der{J}{z\conj} 
			&= \underbrace{\left(\der{J}{s\conj}\right)\conj}_{(\cdot)\in \mathbb{R}} \j\frac{1}{2} + \der{J}{s\conj} \left(-\j\frac{1}{2} \right)\conj \nonumber \\
			&= \gradImag \label{eq:gradImag}
		\end{align}
		yields the gradient.
	\end{rightcolumn}
		
\end{paracol}

\section{DFT/IDFT}
\newcommand{\gradDFT}{}
\newcommand{\gradIDFT}{}

\begin{paracol}{2}
	\setlength{\parskip}{0pt} 
	\setlength{\abovedisplayshortskip}{0pt}
	\setlength{\abovedisplayskip}{0.17in}
	\setlength{\columnseprule}{.4pt}
	\begin{leftcolumn}
		We begin with the definition of the \gls{dft}
		\begin{align}\begin{aligned}
			s_f &= \sum_{n=0}^{N-1} z_n \e^{-\j\frac{2\pi}{N}nf} \\
			&= \dft{z_n}
		\end{aligned}\end{align}
		with $z_n$ as the discrete time signal, $s_f$ as transformed signal in the frequency domain and $n,f\in\left\{0,\dots,N\right\}$.
		Deriving \wrt $z_n$ and $z_n\conj$ 
		\begin{align}
			\der{s_f}{z_n} = \e^{-\j\frac{2\pi}{N}nf} \quad \text{and} \quad
			\der{s_f}{z_n\conj} = 0 \nonumber
		\end{align}
		and inserting in \eqref{eq:chain_rule_tensor} 
		\begin{align}
			\grad{z_n} &= \sum_{f=0}^{N-1} \left(\left(\grad{s_f}\right)\conj 0 + \grad{s_f} \left(\e^{-\j\frac{2\pi}{N}nf}\right)\conj\right) \nonumber \\
			&= \sum_{f=0}^{N-1} \grad{s_f} \e^{\j\frac{2\pi}{N}nf} \nonumber \\
			&= N \cdot \idft{\grad{s_f}} \label{eq:grad:dft}
		\end{align}
		leads to the gradient.
	\end{leftcolumn}
	\begin{rightcolumn}
		We begin with the definition of the \gls{idft}
		\begin{align}\begin{aligned}
			s_n &= \frac{1}{N}\sum_{f=0}^{N-1} z_f \e^{\j\frac{2\pi}{N}nf} \\
			&= \idft{z_f}
		\end{aligned}\end{align}
		with $z_f$ as the discrete signal in the frequency domain, $s_n$ as transformed time signal and $n,f\in\left\{0,\dots,N\right\}$.
		Deriving \wrt $z_f$ and $z_f\conj$ 
		\begin{align}
			\der{s_n}{z_f} = \frac{1}{N} \e^{\j\frac{2\pi}{N}nf} \quad \text{and} \quad
			\der{s_n}{z_f\conj} = 0 \nonumber
		\end{align}
		and inserting in \eqref{eq:chain_rule_tensor} 
		\begin{align}
			\grad{z_f} &= \frac{1}{N}\sum_{n=0}^{N-1} \left(\left(\grad{s_n}\right)\conj 0 + \grad{s_n} \left(\e^{\j\frac{2\pi}{N}nf}\right)\conj\right) \nonumber \\
			&= \frac{1}{N} \sum_{n=0}^{N-1} \grad{s_n} \e^{-\j\frac{2\pi}{N}nf} \nonumber \\
			&= \frac{1}{N} \cdot \dft{\grad{s_n}} \label{eq:grad:idft}
		\end{align}
		leads to the gradient.
	\end{rightcolumn}
\end{paracol}

\section{RDFT/IRDFT}
\newcommand{\gradRDFT}{}
\newcommand{\gradIRDFT}{}
{	

	The size of the \gls{dft} for real input is assumed to be even.
	We begin with the definition of the \gls{rdft}
	\begin{align}\begin{aligned}
		s_f = \sum_{n=0}^{N-1} z_n \e^{-\j\frac{2\pi}{N}nf} = \rdft{z_n}
	\end{aligned}\end{align}
	with $z_n$ as the discrete time signal, $s_f$ as transformed signal in the frequency domain and $n\in\left\{0,\dots,N\right\}, f\in\left\{0,\dots,\frac{N}{2}\right\}$.
	Deriving \wrt $z_n$ and $z_n\conj$ 
	\begin{align}
		\der{s_f}{z_n} = \e^{-\j\frac{2\pi}{N}nf} \quad \text{and} \quad
		\der{s_f}{z_n\conj} = 0 \nonumber
	\end{align}
	and inserting in \eqref{eq:chain_rule_tensor} leads to
	\begin{align}
		\grad{z_n} &= \sum_{f=0}^{N/2} \left(\left(\grad{s_f}\right)\conj 0 + \grad{s_f} \left(\e^{-\j\frac{2\pi}{N}nf}\right)\conj\right) \nonumber \\
		&= \sum_{f=0}^{N/2} \grad{s_f} \e^{\j\frac{2\pi}{N}nf}. \nonumber 
	\end{align}
	Note, that the sum over $f$ is limited by $N/2$.
	Since the input is \rv the gradient has to be \rv, too.
	This leads to
	\begin{align}
	\grad{z_n} &= \Re \left\{\sum_{f=0}^{N/2} \grad{s_f} \e^{\j\frac{2\pi}{N}nf}\right\} \nonumber \\
	&= \sum_{f\in\{0,\frac{N}{2}\}} \Re \left\{\grad{s_f} \right\} \underbrace{\e^{\j\frac{2\pi}{N}nf}}_{\in\mathbb{R}} +
	\sum_{f=1}^{N/2-1} \Re\left\{\grad{s_f} \e^{\j\frac{2\pi}{N}nf}\right\} \nonumber \\
	&= \sum_{f\in\{0,\frac{N}{2}\}} \Re \left\{\grad{s_f} \right\} \underbrace{\e^{\j\frac{2\pi}{N}nf}}_{\in\mathbb{R}} +
	\sum_{f=1}^{N/2-1} \frac{1}{2}\left(\grad{s_f} \e^{\j\frac{2\pi}{N}nf} + \left(\grad{s_f}\right)\conj \e^{-\j\frac{2\pi}{N}nf}\right). \nonumber \\
	\end{align}
	Using the definition of the IRDFT \eqref{eq:irdft} and
	\begin{align}
	\widetilde{\grad{s_f}}	= \begin{cases}
	\Re\left\{\grad{s_f}\right\}, &\text{for } f\in\left\{0, \frac{N}{2}\right\}, \\
	\frac{1}{2} \grad{s_f}, &\text{for } f\in\left[1, \dots, \frac{N}{2}-1\right] \label{eq:grad:rdft1}
	\end{cases}
	\end{align}
	the gradient can be calculated with
	\begin{align}
	\grad{z_n} = N \cdot \irdft{\widetilde{\grad{s_f}}}. \label{eq:grad:rdft2}
	\end{align}

	For the derivation of the \gls{irdft} we start with its definition
	\begin{align}\begin{aligned}
		s_n = \frac{1}{N}\sum_{f=0}^{N/2} z_f \e^{\j\frac{2\pi}{N}nf} 
		 + \frac{1}{N}\sum_{f=1}^{N/2-1} z_f\conj \e^{-\j\frac{2\pi}{N}nf} = \irdft{z_f} \label{eq:irdft}
	\end{aligned}\end{align}
	with $z_f$ as the discrete signal in the frequency domain, $s_n$ as transformed time signal and $n\in\left\{0,\dots,N\right\}, f\in\left\{0,\dots,\frac{N}{2}\right\}$.
	Deriving \wrt $z_f$ and $z_f\conj$ 
	\begin{align}
		\der{s_f}{z_f} = \frac{1}{N} \e^{\j\frac{2\pi}{N}nf} \quad \text{and} \quad
		\der{s_f}{z_f\conj} = \frac{1}{N} \e^{-\j\frac{2\pi}{N}nf} \nonumber
	\end{align}
	and inserting in \eqref{eq:chain_rule_tensor} leads to
	\begin{align}
		\grad{z_0} &= \frac{1}{N} \sum_{n=1}^{N/2-1} \grad{s_n} 1, \nonumber \\
		\grad{z_{N/2}} &= \frac{1}{N} \sum_{n=1}^{N/2-1} \grad{s_n} \e^{\j\pi n}, \nonumber \\
		\grad{z_f} &= \frac{1}{N}\sum_{n=1}^{N/2-1} \left(\left(\grad{s_n}\right)\conj \e^{-\j\frac{2\pi}{N}nf} + \grad{s_n} \left(\e^{\j\frac{2\pi}{N}nf}\right)\conj\right) \nonumber \\
		&= \frac{1}{N}\sum_{n=1}^{N/2-1} 2 \grad{s_n} \e^{-\j\frac{2\pi}{N}nf} \quad \forall z_f, f\in\{1,...,\frac{N}{2}-1\}. \nonumber
	\end{align}
	The last equation holds because the gradient of the objective function \wrt the output of the \gls{irdft} $s_n$ is \rv since the output has to be \rv, too. With the RDFT
	\begin{align}
		\widetilde{\grad{z_f}} = \frac{1}{N} \cdot \rdft{\grad{s_n}} \label{eq:grad:irdft1}
	\end{align}
	the gradient of the \gls{irdft} can be obtained by
	\begin{align}
		\grad{z_f} = \begin{cases}
			\widetilde{\grad{z_f}}, &\text{for } f \in\left\{0, \frac{N}{2}\right\},\\[0.5em]
			2\widetilde{\grad{z_f}}, &\text{for } f \in\left[1, \dots, \frac{N}{2}-1\right]. 
		\end{cases}\label{eq:grad:irdft2}
	\end{align}
}

\section{Vector normalization}
\newcommand{\gradVectorNorm}{}
\label{subsec:vectorNorm}

We define the vector normalization such that the euclidean norm of the output vector is one
\begin{align}
	\vec{s} &= \frac{\vec{z}}{\sqrt{\vec{z}\hermite \vec{z}}}.
\end{align}
The way to determine the gradient with the approach from \cref{sec:der:giles} is more complicated than the approach from \cref{sec:der:drude}. Therefore we start with the elementwise formulation
\begin{align}
	s_n &= \frac{z_n}{\sqrt{\sum_m z_m\conj z_m}} \nonumber,
\end{align}
calculate the derivative \wrt $z_m\conj$
\begin{align}
	\der{s_n}{z_m\conj} 
	&= -\frac{z_n}{\sum_{\check{m}} z_{\check{m}}\conj z_{\check{m}}} \frac{1}{2\sqrt{\sum_{\check{m}} z_{\check{m}}\conj z_{\check{m}}}} z_m 
	= \frac{-s_nz_m}{2\sum_{\check{m}} z_{\check{m}}\conj z_{\check{m}}} \nonumber
\end{align}
and insert it in \eqref{eq:def_der_vec_wrt_vec}
\begin{align}
	\der{\vec{s}}{\vec{z}\conj}
	&= -\frac{\vec{s}\vec{z}\transposed}{2\vec{z}\hermite\vec{z}} \nonumber.
\end{align}
The calculation for $\der{\vec{s}}{\vec{z}}$ is analogue
\begin{align}
	\der{s_n}{z_m} 
	&= -\frac{z_n}{\sum_{\check{m}} z_{\check{m}}\conj z_{\check{m}}} \frac{1}{2\sqrt{\sum_{\check{m}} z_{\check{m}}\conj z_{\check{m}}}} z_m\conj \nonumber
	+ \frac{\delta_{nm}}{\sqrt{\sum_{\check{m}} z_{\check{m}}\conj z_{\check{m}}}},  \nonumber\\
	\der{\vec{s}}{\vec{z}}
	&= \frac{\vec{I}}{\sqrt{\vec{z}\hermite\vec{z}}} - \frac{\vec{s}\vec{z}\hermite}{2\vec{z}\hermite\vec{z}}, \nonumber
\end{align}
where $\delta_{nm}$ is the Kronecker delta
\begin{align}
	\delta_{nm} = 
	\begin{cases}
		1,\quad\mathrm{if}~~m = n,\\
		0,\quad\mathrm{if}~~m \ne n.
	\end{cases}
\end{align}

Inserting this intermediate result in \eqref{eq:chain_rule_vector}
\begin{align}
	\der{J}{\vec{z}\conj} 
	&= \left(\left(\der{J}{\vec{s}\conj}\right)\hermite \der{\vec{s}}{\vec{z}\conj} 
	+ \left(\der{J}{\vec{s}\conj}\right)\transposed \left(\der{\vec{s}}{\vec{z}}\right)\conj\right)\transposed \nonumber\\
	&= \left(-\left(\der{J}{\vec{s}\conj}\right)\hermite \frac{\vec{s}\vec{z}\transposed}{2\vec{z}\hermite\vec{z}} 
	+ \left(\der{J}{\vec{s}\conj}\right)\transposed \left(\frac{\vec{I}}{\sqrt{\vec{z}\hermite\vec{z}}} - \frac{\vec{s}\vec{z}\hermite}{2\vec{z}\hermite\vec{z}}\right)\conj\right)\transposed \nonumber\\
	&= - \frac{\vec{z}\vec{s}\transposed}{2\vec{z}\hermite\vec{z}} \left(\der{J}{\vec{s}\conj}\right)\conj 
	+ \left(\frac{\vec{I}}{\sqrt{\vec{z}\hermite\vec{z}}} - \frac{\vec{z}\conj\vec{s}\transposed}{2\vec{z}\hermite\vec{z}}\right)\conj \der{J}{\vec{s}\conj} \nonumber\\
	&= \frac{1}{\sqrt{\vec{z}\hermite\vec{z}}} \der{J}{\vec{s}\conj} 
	- \frac{\vec{z}\vec{s}\transposed}{2\vec{z}\hermite\vec{z}} \left(\der{J}{\vec{s}\conj}\right)\conj 
	- \frac{\vec{z}\vec{s}\hermite}{2\vec{z}\hermite\vec{z}} \der{J}{\vec{s}\conj} \nonumber\\
	&= \frac{1}{\sqrt{\vec{z}\hermite\vec{z}}} \der{J}{\vec{s}\conj} 
	- \frac{\vec{z}}{2\vec{z}\hermite\vec{z}} \left(\vec{s}\transposed\left(\der{J}{\vec{s}\conj}\right)\conj +\vec{s}\hermite\der{J}{\vec{s}\conj} \right) \nonumber\\
	&= \frac{1}{\sqrt{\vec{z}\hermite\vec{z}}} \der{J}{\vec{s}\conj} 
	- \frac{\vec{z}}{\vec{z}\hermite\vec{z}} \Re\left\{\vec{s}\hermite\der{J}{\vec{s}\conj} \right\} \nonumber\\
	&= \frac{1}{\sqrt{\vec{z}\hermite\vec{z}}} \der{J}{\vec{s}\conj} 
	- \frac{\vec{z}}{\sqrt{\vec{z}\hermite\vec{z}}^3} \Re\left\{\vec{z}\hermite\der{J}{\vec{s}\conj} \right\} \nonumber\\
	&= \frac{\der{J}{\vec{s}\conj}-\frac{\vec{z}}{\vec{z}\hermite\vec{z}}\Re\{\vec{z}\hermite\der{J}{\vec{s}\conj}\}}{\sqrt{\vec{z}\hermite\vec{z}}} \label{eq:grad:vecNormalisation}
\end{align}
yields the gradient.

\section{Matrix multiplication}
\label{sec:matMul}
\newcommand{\gradMatMul}{}

For the matrix multiplication
\begin{align}
	\vec{C} = \vec{g}(\vec{A}, \vec{B}) = \vec{A}\vec{B}
\end{align}
we want to show two ways for the calculation of the gradient. In the left column we show the way like in the previous sections (\cref{sec:der:drude}) and use the elementwise notation
\begin{align}
	c_{nm} &= g(a_{n1}, \dots, a_{nK}, b_{1m}, \dots, b_{Km}) \nonumber \\
	&= \sum_k a_{nk}b_{km} \nonumber.
\end{align}
Note, that here the variable names have changed. In the previous sections the inputs/outputs are $z_1$, $z_2$ and $s$ and now they are $\vec{A}$, $\vec{B}$ and $\vec{C}$. This change is done, because the algorithm in \cref{sec:der:giles} requires a pseudo input variable, which we call $\zeta \in \{z, z\conj\}$. This algorithm is used in the right column.
\pagebreak

\columnratio{0.3}
\begin{paracol}{2}
	\setlength{\parskip}{0pt} 
	\setlength{\columnseprule}{.4pt}
	\begin{leftcolumn}
		We start with elementwise notation
		\begin{align}
			c_{n,m} &= \sum_k a_{nk}b_{k,m}, \nonumber
		\end{align}	
		calculate all \\ necessary derivatives
		\begin{align}
		\der{c_{n,m}}{a_{n,k}} &= b_{k,m}, \nonumber\\
		\der{c_{n,m}}{a_{n,k}\conj} &= 0, \nonumber\\
		\der{c_{n,m}}{b_{k,m}} &= a_{n,k}, \nonumber\\
		\der{c_{n,m}}{b_{k,m}\conj} &= 0. \nonumber
		\end{align}	
		and inserting them in \eqref{eq:chain_rule_tensor}
		\begin{align}
		\grad{a_{n,k}} &= \left(\grad{c_{n,m}}\right)\conj 0 \nonumber \\
		&+ \grad{c_{n,m}} b_{k,m}\conj \nonumber \\
		\grad{\vec{A}} &= \grad{\vec{C}} \vec{B}\hermite \\
		\grad{b_{k,m}} &= \left(\grad{c_{n,m}}\right)\conj 0 \nonumber\\
		&+ \grad{c_{n,m}} a_{nk}\conj \nonumber \\
		\grad{\vec{B}} &= \vec{A}\hermite\grad{\vec{C}}
		\end{align}
		yields the result.
	\end{leftcolumn}
	\begin{rightcolumn}
		Start with calculating the derivative\footnote{In \cref{appendix:proof:complex_valued_chain_rule} it is shown, that it is not necessary to distinguish between $z$ and $z\conj$ in this case.} of both sides of $\vec{C}=\vec{A}\vec{B}$ \wrt $\zeta \in \{z, z\conj\}$. 
		In \autoref{subsec:matrixproduct} it is shown that the product rule for real matrices holds also in the complex case \wrt $z$ and $z\conj$, which gives
		\begin{align}
		\derz{\vec{C}} = \derz{\vec{A}}\vec{B} + \vec{A}\derz{\vec{B}}. \nonumber
		\end{align}
		Insert this in \eqref{eq:giles:trace_c}
		\begin{align}
		\grad{z} &= \trace{\left(\grad{\vec{C}}\right)\hermite \der{\vec{C}}{z\conj}} + \dots \nonumber\\
		&= \trace{\left(\grad{\vec{C}}\right)\hermite \left(\der{\vec{A}}{z\conj}\vec{B} + \vec{A}\der{\vec{B}}{z\conj}\right)} + \dots \nonumber\\
		&= \trace{\left(\grad{\vec{C}}\vec{B}\hermite\right)\hermite \der{\vec{A}}{z\conj}
			+ \left(\vec{A}\hermite\grad{\vec{C}}\right)\hermite \der{\vec{B}}{z\conj}} + \dots \nonumber
		\end{align}
		The dots represent the part of \eqref{eq:giles:trace_c}, which is redundant in this example.
		Comparing the coefficients with \eqref{eq:giles:trace_a_b}
		\begin{align}
		\grad{\vec{A}} &= \grad{\vec{C}}\vec{B}\hermite \label{eq:grad:matMul1} \\
		\grad{\vec{B}} &= \vec{A}\hermite\grad{\vec{C}} \label{eq:grad:matMul2}
		\end{align}
		yields the results.
	\end{rightcolumn}
\end{paracol}
As expected both ways lead to the same solution.

\section{Matrix inverse}
\newcommand{\gradMatInv}{}
\label{subsec:matrixInv}

For the matrix inverse
\begin{align}
	\vec{C} = \vec{g}(\vec{A}) = \vec{A}\inv
\end{align}
we assume the inverse exists. An implicit formulation is
\begin{align}
	\vec{A}\vec{C} = \vec{I} \nonumber
\end{align}
and the derivation \wrt $\zeta \in \{z, z\conj\}$ leads to
\begin{align}
	\derz{\vec{A}}\vec{C} + \vec{A}\derz{\vec{C}} = \vec{0}, \nonumber\\
	\derz{\vec{C}} = -\vec{A}\inv\derz{\vec{A}}\vec{C}.\nonumber
\end{align}
Inserting these equations in \eqref{eq:giles:trace_c}
\begin{align}
	\grad{z} &= \trace{\left(\grad{\vec{C}}\right)\hermite \der{\vec{C}}{z\conj}} + \dots \nonumber\\
	&= \trace{\left(\grad{\vec{C}}\right)\hermite  \left(-\vec{A}\inv\der{\vec{A}}{z\conj}\vec{C}\right) } + \dots \nonumber\\
	&= \trace{-\vec{C}\left(\grad{\vec{C}}\right)\hermite  \vec{A}\inv\der{\vec{A}}{z\conj} } + \dots \nonumber\\
	&= \trace{\left(-\vec{A}\hermiteinv\grad{\vec{C}}\vec{C}\hermite\right)\hermite  \der{\vec{A}}{z\conj} } + \dots \nonumber\\
	&= \trace{\left(-\vec{C}\hermite\grad{\vec{C}}\vec{C}\hermite\right)\hermite  \der{\vec{A}}{z\conj} } + \dots\nonumber
\end{align}
and comparing the coefficients with \eqref{eq:giles:trace_a_b}
\begin{align}
	\grad{\vec{A}} &= -\vec{C}\hermite\grad{\vec{C}}\vec{C}\hermite \label{eq:grad:matInv}
\end{align}
yields the gradient.

\section{Matrix inverse product}
\label{subsec:matrixProdInv}
\newcommand{\gradMatMulInv}{}

The matrix inverse product
\begin{align}
	\vec{C} = \vec{g}(\vec{A}, \vec{B}) = \vec{A}\inv\vec{B}
\end{align}
is also called \texttt{solve} in many numerical implementations, because it solves the equation
\begin{align}
	\vec{A}\vec{C} = \vec{B}. \nonumber
\end{align}
The derivative of this equation \wrt to $\zeta \in \{z, z\conj\}$ is
\begin{align}
	\derz{\vec{A}}\vec{C} + \vec{A}\derz{\vec{C}} = \derz{\vec{B}}, \nonumber\\
	\derz{\vec{C}} = \vec{A}\inv\derz{B} - \vec{A}\inv\derz{\vec{A}}\vec{C}. \nonumber
\end{align}
Now the same principle as in the previous section can be used.
Inserting in \eqref{eq:giles:trace_c}
\begin{align}
	\grad{z} &= \trace{\left(\grad{\vec{C}}\right)\hermite \der{\vec{C}}{z\conj}} + \dots \nonumber \\
	&= \trace{\left(\grad{\vec{C}}\right)\hermite  \vec{A}\inv\left(\der{\vec{B}}{z\conj} - \der{\vec{A}}{z\conj}\vec{C}\right) } + \dots \nonumber\\
	&= \trace{\left(\vec{A}\hermiteinv\grad{\vec{C}}\right)\hermite \der{\vec{B}}{z\conj} + \left(-\vec{A}\hermiteinv\grad{\vec{C}}\vec{C}\hermite\right)\hermite\der{\vec{A}}{z\conj} } + \dots \nonumber
\end{align}
and comparing with \eqref{eq:giles:trace_a_b}
\begin{align}
	\grad{\vec{B}} &= \vec{A}\hermiteinv\grad{\vec{C}} \label{eq:grad:invMatMul11} \\
	\grad{\vec{A}} &= -\vec{A}\hermiteinv\grad{\vec{C}}\vec{C}\hermite = -\grad{\vec{B}}\vec{C}\hermite \label{eq:grad:invMatMul12}
\end{align}
yields the result.

For the second matrix inverse product
\begin{align}
	\vec{C} = \vec{g}(\vec{A}, \vec{B}) = \vec{A}\vec{B}\inv
\end{align}
the calculation is analogue.
Start with implicit formulation 
\begin{align}
	\vec{C}\vec{B} = \vec{A}, \nonumber
\end{align}
calculate the derivative
\begin{align}
	\derz{\vec{C}}\vec{B} + \vec{C}\derz{\vec{B}} = \derz{\vec{A}}, \nonumber\\
	\derz{\vec{C}} = \derz{A}\vec{B}\inv - \vec{C}\derz{\vec{B}}\vec{B}\inv,\nonumber
\end{align}
insert it in \eqref{eq:giles:trace_c}
\begin{align}
	\grad{z} &= \trace{\left(\grad{\vec{C}}\right)\hermite \der{\vec{C}}{z\conj}} + \dots \nonumber\\
	&= \trace{\left(\grad{\vec{C}}\right)\hermite  \left(\der{A}{z\conj} - \vec{C}\der{\vec{B}}{z\conj}\right)\vec{B}\inv } + \dots \nonumber\\
	&= \trace{\left(\grad{\vec{C}}\vec{B}\hermiteinv\right)\hermite \der{\vec{A}}{z\conj} + \left(-\vec{C}\hermite\grad{\vec{C}}\vec{B}\hermiteinv\right)\hermite\der{\vec{B}}{z\conj} } + \dots \nonumber
\end{align}
and compare with \eqref{eq:giles:trace_a_b}
\begin{align}
	\grad{\vec{A}} &= \grad{\vec{C}}\vec{B}\hermiteinv, \label{eq:grad:invMatMul21}\\
	\grad{\vec{B}} &= -\vec{C}\hermite\grad{\vec{C}}\vec{B}\hermiteinv = -\vec{C}\hermite\grad{\vec{A}} \label{eq:grad:invMatMul22}
\end{align}
which yields the result again.

\section{Cholesky decomposition}
\newcommand{\gradCholesky}{}

The derivation for the \rv Cholesky decomposition is described in \cite[Chapter 3]{Murray2016CholeskyAD} and the following \cv calculation is similar to it.
The Cholesky decomposition is denoted as 
\begin{align}
	\vec{L} = \vec{g}(\vec{A}) = \mathrm{cholesky}\left(\vec{A}\right)
\end{align}
where $ \vec{A} = \vec{L}\vec{L}\hermite$ and $\vec{L}$ is a lower triangular matrix 
(in the following $L \in \llblacktriangle$ and $U \in \urblacktriangle$ indicate a lower and upper triangular matrix, respectively). 
Also we assume, that $\vec{A}$ is a hermitian matrix ($\vec{A}=\vec{A}\hermite$). 
Therefore the operation $\vec{A} \stackrel{!}{=} \frac{1}{2}\left(\vec{B}+\vec{B}\hermite\right)$ where $\vec{B} = \vec{A}$ holds. 
The gradient of $\vec{B}$ is then $\grad{\vec{B}} = \frac{1}{2}\left(\grad{\vec{A}}+\left(\grad{\vec{A}}\right)\hermite\right)$ (see \cref{sec:block:add,sec:block:conj}).

Some properties of triangular matrices and the Cholesky decomposition:
\begin{enumerate}[label=\textbf{Prop.\arabic*)},ref={Prop.\arabic*},align=parleft,leftmargin=*]
	\newcounter{NameOfTheNewCounter}
	\item
		If $\vec{L} \in \llblacktriangle$ then $\vec{L}\inv \in \llblacktriangle$.
	\item
		If $\vec{L}_1 \in \llblacktriangle$ and $ \vec{L}_2 \in \llblacktriangle$ then $\vec{L}_1\vec{L}_2 \in \llblacktriangle$.
	\item
		If $\vec{L} \in \llblacktriangle$ then $\left[\vec{L}\inv\right]_{i,i} = \left(L_{i,i}\right)\inv$ ($\left[\cdot\right]_{i,j}$ takes element from row $i$ and column $j$).
	\item
		If $\vec{L}_1, \vec{L}_2, \vec{L}_3 \in \llblacktriangle$ and $\vec{L}_3 = \vec{L}_1 \vec{L}_2$ then $\left[\vec{L}_3\right]_{i,i} = \left[\vec{L}_1\right]_{i,i} \left[\vec{L}_2\right]_{i,i}$ .
	\item
		diagonal entries of $\vec{L} = \vec{g}(\vec{A})$ are real and positive (from previous follows: diagonal entries of $\vec{L}\inv$ are also real and positive).
	\item
		If $\vec{L}_1, \vec{L}_2 \in \llblacktriangle$, $A_2 = \vec{L}_2\vec{L}_2\hermite$ and $\vec{A}_3 = \vec{L}_1\vec{A}_2\vec{L}_1\hermite$ then $\vec{A}_3 = \vec{L}_1\vec{L}_2\left(\vec{L}_1\vec{L}_2\right)\hermite $ is a hermitian matrix ($\vec{A}_3 = \vec{A}_3\hermite$).
	\item
		If $\vec{L} \in \llblacktriangle$, $L_{i,i} \in \mathbb{R}$, $\vec{C} = \vec{C}\hermite$ and $\vec{C} = \vec{L} + \vec{L}\hermite$ then $\left(\lltriangle[very thick] - \frac{1}{2}\vec{I}\right) \hadamard \vec{C} = \vec{L}$, where \lltriangle[very thick] is a lower triangular matrix filled with ones on the diagonal and below and $\hadamard$ denotes the Hadamard product (also known as elementwise multiplication). \label{item:test}
\end{enumerate}

We start with the equation $\vec{L} = \vec{g}(\vec{A})$ and calculate the derivative \wrt $\zeta$
\begin{align}
	\derz{\vec{A}} &= \derz{\vec{L}}\vec{L\hermite} + \vec{L}\derz{\vec{L\hermite}}. \nonumber
\end{align}
After reformulation we obtain:
\begin{align}
	\vec{L}\inv\derz{\vec{A}}\vec{L}\hermiteinv &= \vec{L}\inv\derz{\vec{L}} + \derz{\vec{L}\hermite}\vec{L}\hermiteinv \nonumber, \\
	\vec{L}\inv\derz{\vec{A}}\vec{L}\hermiteinv &= \underbrace{\vec{L}\inv\derz{\vec{L}}}_{(\cdot) \in \llblacktriangle} + \underbrace{\left(\vec{L}\inv\derz{\vec{L}}\right)\hermite}_{(\cdot) \in \urblacktriangle}. \nonumber
\end{align}
Assuming that the imaginary part of the diagonal of $\vec{L}\inv\derz{\vec{L}}$ is zero\footnote{
	The assumption that the imaginary part of the diagonal of $\vec{L}\inv\derz{\vec{L}}$ has to be zero comes from \ref{item:test}. 
	We have no analytical solution, that this is the case.
	But the final solution for the gradient survives a numerical test. 
	We know: If and only if the diagonal from $\derz{\vec{L}}$ is real valued this is the case.}, 
we can reformulate the equation:
\begin{align}
	\left(\lltriangle[very thick] - \frac{1}{2}\vec{I}\right) \hadamard \left(\vec{L}\inv\derz{\vec{A}}\vec{L}\hermiteinv\right) &= \vec{L}\inv\derz{\vec{L}}, \nonumber\\
	\derz{\vec{L}} &= \vec{L}\left(\left(\lltriangle[very thick] - \frac{1}{2}\vec{I}\right) \hadamard \left(\vec{L}\inv\derz{\vec{A}}\vec{L}\hermiteinv\right)\right). \nonumber
\end{align}
Inserting in \eqref{eq:giles:trace_c}
\begin{align}
	\grad{z} &= \trace{\left(\grad{\vec{L}}\right)\hermite\der{\vec{L}}{z\conj}} + \dots \nonumber \\
	&= \trace{\left(\grad{\vec{L}}\right)\hermite\vec{L}\left(\left(\lltriangle[very thick] - \frac{1}{2}\vec{I}\right) \hadamard \left(\vec{L}\inv\derz{\vec{A}}\vec{L}\hermiteinv\right)\right)} + \dots \nonumber
\end{align}
and using
$\trace{\vec{A}\left(\vec{B} \hadamard \vec{C}\right)} = \trace{\left(\vec{A} \hadamard \vec{B}\transposed\right)\vec{C}}$ \cite[Theorem 2.6]{million2007hadamard}
results in:
\begin{align}
	\grad{z} 
	&= \trace{\left(\left(\left(\grad{\vec{L}}\right)\hermite\vec{L}\right)\hadamard \left(\lltriangle[very thick] - \frac{1}{2}\vec{I}\right)\transposed \right) \vec{L}\inv\derz{\vec{A}}\vec{L}\hermiteinv} + \dots \nonumber\\
	&= \trace{\vec{L}\hermiteinv\left(\left(\left(\grad{\vec{L}}\right)\hermite\vec{L}\right)\hadamard \left(\lltriangle[very thick] - \frac{1}{2}\vec{I}\right)\transposed \right) \vec{L}\inv\derz{\vec{A}}} + \dots \nonumber\\
	&= \trace{\left(\vec{L}\hermiteinv\left(\left(\vec{L}\hermite\grad{\vec{L}}\right)\hadamard \left(\lltriangle[very thick] - \frac{1}{2}\vec{I}\right) \right)\vec{L}\inv\right)\hermite \derz{\vec{A}}} + \dots . \nonumber
\end{align}
With the comparison \eqref{eq:giles:trace_a_b} we receive
\begin{align}
	\grad{\vec{A}} &= \vec{L}\hermiteinv\left(\left(\vec{L}\hermite\grad{\vec{L}}\right)\hadamard \left(\lltriangle[very thick] - \frac{1}{2}\vec{I}\right) \right)\vec{L}\inv.
	\label{eq:chol:grad}
\end{align}
Taking into account, that the input for the Cholesky decomposition must be a hermitian matrix, 
it is intuitive to ensure that the gradient is also a hermitian matrix 
(Note: in general $\grad{\vec{A}}$ is not a hermitian matrix)
\begin{align}
	\grad{\vec{B}} &= \frac{1}{2}\left(\grad{\vec{A}} + \left(\grad{\vec{A}}\right)\hermite\right).
	\label{eq:chol:gradStable}
\end{align}
This force hermitian operation is debatable, 
because simulations have shown, that for a calculation, 
which produces a hermitian matrix, 
it makes no difference for the gradient if we take $\grad{\vec{A}}$ or $\grad{\vec{B}}$. 
On the contrary, there is an analytical interpretation, 
where $\vec{A}$ is a parameter, 
which needs an update $\vec{A}^\mathrm{new}=\vec{A}^\mathrm{old}-\mu \grad{\vec{A}}$. 
Here the gradient must ensure that $\vec{A}$ remains hermitian. 
To ensure this, the gradient has to be hermitian. 
Therefore we use $\grad{\vec{B}}$.



\section{Eigendecomposition}
To perform the Beamforming operation in a \gls{NN} we need the derivate of the eigendecomposition with complex-valued eigenvectors. 
This derivation is similar to the real valued version in \cite{Giles2008}. 
We start with the eigenvalue equation
\begin{align}
	\vecG{\Phi} \vec{W} = \vec{W} \vecG{\Lambda}
\end{align}
where $\vecG{\Phi}$, $\vec{W}$ and $\vecG{\Lambda}$ are a quadratic matrix, the eigenvector matrix and the diagonal eigenvalue matrix, respectively. 
At first we examine the forward mode \gls{ad}.
With the matrix product rule the derivation can be written as
\begin{align}
	\derz{\vecG{\Phi}} \vec{W} + \vecG{\Phi} \derz{\vec{W}} &= \derz{\vec{W}} \vecG{\Lambda} + \vec{W} \derz{\vecG{\Lambda}}, \label{eq:eig:ausgang}
\end{align}
where $\zeta \in \{z,z\conj\}$. 
Since the matrix of eigenvalues is a diagonal matrix we use the rule (proof in \cref{appendix:eigen:CLambdaMinusLamdaC})
\begin{align}
	\vec{C} \vecG{\Lambda} - \vecG{\Lambda}\vec{C} = \vec{E} \circ \vec{C},
	\label{eq:eig:tmp1}
\end{align}
where $(A\circ B)_{ij} = A_{ij}B_{ij}$ denotes the Hadamard product and $E_{ij} = \lambda_j - \lambda_i$ are the entries of $\vec{E}$, which is a matrix of the differences of the eigenvalues \cite{Giles2008}. 

With the identity $\vecG{\Phi} = \vec{W}\vecG{\Lambda}\vec{W}\inv$, \eqref{eq:eig:ausgang} premultiplied with $\vec{W}\inv$
\begin{align}\begin{aligned}
	\vec{W}\inv \derz{\vecG{\Phi}} \vec{W} - \derz{\vecG{\Lambda}} &=  \vec{W}\inv\derz{\vec{W}} \vecG{\Lambda} - \vecG{\Lambda}\vec{W}\inv\derz{\vec{W}}
\end{aligned}\end{align}
 and \eqref{eq:eig:tmp1} we get
\begin{align}
	\vec{W}\inv \derz{\vecG{\Phi}} \vec{W} - \derz{\vecG{\Lambda}} &=  \vec{E} \circ ( \vec{W}\inv  \derz{\vec{W}}). \label{eq:eig:tmp2}
\end{align}

After applying the Hadamard product with $\vec{I}$ to both sides of \eqref{eq:eig:tmp2} the forward mode gradient of the eigenvalues is calculated via
\begin{align}
	\derz{\vecG{\Lambda}} = \vec{I} \circ (\vec{W}\inv \derz{\vecG{\Phi}} \vec{W}),
	\label{eq:eig:derEW}
\end{align}
because $\vec{E}\hadamard\vec{I} = 0$.
Elementwise multiplication of \eqref{eq:eig:tmp2} with $(\vec{F} \circ)$, which is the Hadamard inverse of $\vec{E}$ except on the diagonal where it is zero, leads to
\begin{align}
	\underbrace{\vec{F} \circ \vec{E}}_{\vec{1} - \vec{I}} \circ ( \vec{W}\inv  \derz{\vec{W}}) &= \vec{F} \circ \left(\vec{W}\inv \derz{\vecG{\Phi}} \vec{W}\right) - \underbrace{\vec{F} \circ \derz{\vecG{\Lambda}}}_{0}. \label{eq:eig:tmp3}
\end{align}
The diagonal elements of $\vec{W}\inv \derz{\vec{W}}$ have to be zero, because a change of $\zeta$ can not influence the amplitude of the eigenvectors \cite{Giles2008}. 
Therefore $\vec{I} \circ ( \vec{W}\inv  \derz{\vec{W}})$ is zero and 
\begin{align}
	\vec{W}\inv \derz{\vec{W}} &= \vec{F} \circ \left(\vec{W}\inv \derz{\vecG{\Phi}} \vec{W}\right), \\
	\derz{\vec{W}} &= \vec{W} \left[\vec{F} \circ \left(\vec{W}\inv \derz{\vecG{\Phi}} \vec{W}\right)\right]
	\label{eq:eig:derEV}
\end{align}
denotes the forward mode gradient of the eigenvectors.

To calculate the backward mode gradient of the eigendecomposition, which refers to the derivative of the objective function \wrt the \gls{PSD} matrix, we have to go all the way back through the forward mode calculation. 
Since \eqref{eq:complChainRule} can be reformulated to 
\begin{align}\begin{aligned}
	\der{J}{z\conj} &= \trace{\left(\der{J}{\vecG{\Phi}\conj}\right)\hermite \der{\vecG{\Phi}}{z\conj}} +  \trace{\left(\der{J}{\vecG{\Phi}\conj}\right)\transposed\left(\der{\vecG{\Phi}}{z}\right)\conj} \\
	&= \trace{\left(\der{J}{\vecG{\Phi}\conj}\right)\hermite \der{\vecG{\Phi}}{z\conj}} +  \trace{\left(\der{J}{\vecG{\Phi}\conj}\right)\hermite\der{\vecG{\Phi}}{z}}\conj.
\end{aligned}\end{align}
it is sufficient to focus on the first trace
\begin{align}
	\der{J}{z\conj} = \trace{\left(\der{J}{\vecG{\Phi}\conj}\right)\hermite \der{\vecG{\Phi}}{z\conj}} + \dots
	\label{eq:tracePhi}
\end{align}
because the differentiation \wrt $z$ leads to the same solution.
The dots represent the neglected differentiation \wrt $z$.
The chain rule \eqref{eq:complChainRule} also holds for two intermediate matrices:
\begin{align}
	\der{J}{z\conj} = \trace{\left(\der{J}{\vec{W}\conj}\right)\hermite \der{\vec{W}}{z\conj} + \left(\der{J}{\vecG{\Lambda}\conj}\right)\hermite \der{\vecG{\Lambda}}{z\conj}} + \dots
	\label{eq:traceLambdaW}
\end{align}
Now we can compare the coefficients of \eqref{eq:tracePhi} with \eqref{eq:traceLambdaW} to calculate $\der{J}{\vecG{\Phi}\conj}$.
Inserting the forward mode gradients \eqref{eq:eig:derEW} and \eqref{eq:eig:derEV}
\begin{align}\begin{aligned}
	\der{J}{z\conj} &= \trace{\left(\der{J}{\vec{W}\conj}\right)\hermite \vec{W} \left[\vec{F} \circ \left(\vec{W}\inv \der{\vecG{\Phi}}{z\conj} \vec{W}\right)\right]} \\
	&+ \trace{\left(\der{J}{\vecG{\Lambda}\conj}\right)\hermite \vec{I} \circ (\vec{W}\inv \der{\vecG{\Phi}}{z\conj} \vec{W})}  + \dots 
\end{aligned}\end{align}
and rearranging the matrices with 
\begin{align}
\Tr\left(AB\right) &= \Tr\left(BA\right), \\
\trace{A(B\circ C)} &= \trace{(A\circ B\transposed)C} \quad \text{\cite[Theorem 2.6]{million2007hadamard}},
\end{align}
leads to 
\begin{align}
	\der{J}{z\conj} &= \trace{\vec{W}\left[\left(\der{J}{\vec{W}\conj}\right)\hermite \vec{W} \circ \vec{F}\transposed+ \left(\der{J}{\vecG{\Lambda}\conj}\right)\hermite\right]\vec{W}\inv\der{\vecG{\Phi}}{z\conj} }  
	 + \dots
\end{align}
The derivation of the objective function \wrt $\vecG{\Phi}\conj$ can now be obtained via comparing the coefficients. This leads to the backward mode gradient of the eigendecomposition
\begin{align}
	\der{J}{\vecG{\Phi}\conj} &= \left(\vec{W}\left[\left(\der{J}{\vec{W}\conj}\right)\hermite \vec{W} \circ \vec{F}\transposed+ \left(\der{J}{\vecG{\Lambda}\conj}\right)\hermite\right]\vec{W}\inv\right)\hermite \nonumber \\
	&= \vec{W}\hermiteinv \left[\der{J}{\vecG{\Lambda}\conj} + \vec{F}\conj \circ \vec{W}\hermite \der{J}{\vec{W}\conj} \right] \vec{W}\hermite \label{eq:eig:backwardgradient},
\end{align}
where $(\cdot)\hermiteinv$ indicates the inverse of the conjugate transpose.
In \cite{Giles2008} the real-valued version of this equation can be found.

\subsection{Extension of the gradient}
The magnitude of the eigenvectors is one degree of freedom of the eigendecomposition. 
It is set to one in many numerical implementations.
If the following calculations depend on the magnitude of the eigenvectors and $\vecG{\Phi}$ is not a hermitian matrix this has to be considered in the backward gradient.
Here we use the normalization for vectors (\autoref{subsec:vectorNorm})
\begin{align}
\vec{a} &= \frac{\vec{w}}{\sqrt{\vec{w}\hermite \vec{w}}}, \\
\der{J}{\vec{w}\conj} \label{eq:eig:norm} &= \frac{\der{J}{\vec{a}\conj}-\vec{w}\frac{\Re\{\vec{w}\hermite\der{J}{\vec{a}\conj}\}}{\vec{w}\hermite\vec{w}}}{\sqrt{\vec{w}\hermite\vec{w}}}.
\end{align}

If $\vec{w}$ has already the magnitude one ($\vec{a} = \vec{w}$), this simplifies the equation to
\begin{align}
\der{J}{\vec{w}\conj} &= \der{J}{\vec{a}\conj}-\vec{w}\Re\{\vec{w}\hermite\der{J}{\vec{a}\conj}\}.
\end{align}
The extension to all eigenvectors $\vec{W} = [\vec{w}_1, \dots, \vec{w}_N]$ leads to
\begin{align}
\der{J}{\vec{W}\conj} &= \der{J}{\vec{A}\conj}-\vec{W}\left(\Re\{\vec{W}\hermite\der{J}{\vec{A}\conj}\}\circ\vec{I}\right)
\end{align}
and the complete gradient of the eigendecomposition is
\begin{align}\begin{aligned}
\der{J}{\vecG{\Phi}\conj} &= \vec{W}\hermiteinv \left(\der{J}{\vecG{\Lambda}\conj} + \vec{F}\conj \circ \left(\vec{W}\hermite \der{J}{\vec{A}\conj} \right)\right) \vec{W}\hermite \\
&-\vec{W}\hermiteinv \left(\vec{F}\conj \circ \vec{W}\hermite\vec{W}\left(\Re\{\vec{W}\hermite\der{J}{\vec{A}\conj}\}\circ\vec{I}\right)\right) \vec{W}\hermite. \label{eq:eig:prefinal}
\end{aligned}\end{align}
Since \eqref{eq:eig:norm} is not actually calculated during the forward step, we can relabel $\grad{\vec{A}}$ back into $\grad{\vec{W}}$
\begin{align}\begin{aligned}
\der{J}{\vecG{\Phi}\conj} &= \vec{W}\hermiteinv \left(\der{J}{\vecG{\Lambda}\conj} + \vec{F}\conj \circ \left(\vec{W}\hermite \der{J}{\vec{W}\conj} \right)\right) \vec{W}\hermite \\
&-\vec{W}\hermiteinv \left(\vec{F}\conj \circ \vec{W}\hermite\vec{W}\left(\Re\{\vec{W}\hermite\der{J}{\vec{W}\conj}\}\circ\vec{I}\right)\right) \vec{W}\hermite. \label{eq:eig:final}
\end{aligned}\end{align}

This extension is not always necessary.
When $\vecG{\Phi}$ is a hermitian matrix, then the eigenvector matrix $\vec{W}$ is an orthogonal matrix ($\vec{W}\inv = \vec{W}\hermite$).
With this constraint 
\begin{align}
	\vec{I} = \vec{W}\hermite\vec{W}
\end{align}
holds and since $\vec{F}$ is an off diagonal matrix, then
\begin{align}
	\vec{0} = \vec{F}\conj \hadamard \vec{I} \hadamard \left(\cdot\right),
\end{align}
and the extension is zero.

A second case is, when the following calculations do not depend on the magnitude of the eigenvectors.
Then the dot product between an eigenvector and its gradient
\begin{align}
\Re\{\vec{w}_n\hermite\der{J}{\vec{w}_n\conj}\}
\end{align}
is zero, because the projection of the gradient in the direction of the eigenvector is zero, i.e. the angle between them is $\SI{90}{\degree}$.
Since this values are the diagonal elements of 
\begin{align}
	\Re\{\vec{W}\hermite\der{J}{\vec{W}\conj}\}\circ\vec{I},
\end{align}
the extension is zero. 

\subsection{Verification of the gradient}
The eigendecomposition is invertible.
Therefore the calculation for
\begin{align}
	\vec{\Lambda}, \vec{W} &= \mathrm{eig}(\vecG{\phi}), \\
	\tilde{\vecG{\phi}} &= \vec{W} \vec{\Lambda} \vec{W}\inv,
\end{align}
with an arbitrary matrix $\vecG{\phi}$ and backward gradient $\der{J}{\tilde{\vecG{\phi}}\conj}$ should not affect the gradient. The proof of this can be found in \cref{appendix:proof:eig_identity}.
One conclusion of this calculation is, that the extension of the gradient for this example is redundant, because it is zero.
Also in a numerical verification this identity holds.

For the verification of the necessity of the extension the following numerical proof is done, by calculating the numeric analytical and approximated numeric (\cref{sec:numeric:derivative}) gradient of 
\begin{align}
\vec{\Lambda}, [\vec{w}_1, \dots, \vec{w}_N] &= \mathrm{eig}(\vecG{\phi}), \\
\tilde{\vec{w}}_n &= \vec{w}_n \e^{-\j \angle w_{1n}} \quad \text{for}~n \in\{1, \dots, N\},
\end{align}
where $\angle$ extracts the angle. 
The second equation removes the degree of freedom in an absolute phase of the eigenvectors. 
This is necessary, because the approximated numeric gradient requires a similar output for similar inputs.
The simulations have shown, that for this test the extended gradient is necessary to achieve the correct gradient.

The third verification is again a numerical test of the gradient, but this time in the context of the \gls{gev} beamformer.
The \gls{gev} beamformer with a unique output can be written as follows:
\begin{align}
	\vec{\Lambda}, [\vec{w}_1, \dots, \vec{w}_N] &= \mathrm{gev}\left(\sum_t M^{(X)}_t \vec{Y}_t\vec{Y}_t\hermite, \sum_t M^{(N)}_t \vec{Y}_t\vec{Y}_t\hermite\right), \\
	\vec{w} &= \vec{w}_n \e^{-\j \angle w_{1n}},
\end{align}
where $n=\argmax{n}{\lambda_n}$ and $\lambda_n$ are the diagonal entries of $\vec{\Lambda}$ and the function gev solves the generalized eigenvalue problem.
An alternative calculation of the beamformer is over spatial whitening (\cref{appendix:GEV_spatial_whitening})
\begin{align}
	\phinn^{\frac{1}{2}} &= \text{cholesky}\left(\sum_t M^{(N)}_t \vec{Y}_t\vec{Y}_t\hermite\right), \\
	\tilde{\vec{Y}}_t &= \phinn^{-\frac{1}{2}} \vec{Y}_t, \\
	\vec{\Lambda}, [\tilde{\vec{w}}_1, \dots, \tilde{\vec{w}}_N] &= \mathrm{eig}\left(\sum_t M^{(X)}_t \tilde{\vec{Y}}_t^{(X)}\tilde{\vec{Y}}_t\hermite\right), \\
	\vec{w}_n &= \left(\phinn^{-\frac{1}{2}}\right)\hermite \tilde{\vec{w}}_n, \\
	\vec{w} &= \frac{\vec{w}_n}{\left|\vec{w}_n\right|} \e^{-\j \angle w_{1n}},
\end{align}
where only the gradient for the eigenvalue problem for hermitian matrices is required.
As it should be, both beamformers are equal and also the gradients for $M^{(X)}_t$, $M^{(N)}_t$ and $\vec{Y}_t$ are equal.
Of cause the gradient also survives a numerical verification.

\pagebreak


\section{Overview}
\label{sec:blockSummary}
In the end of this chapter all gradients are summarized in three tables (\cref{tab:grad:scalar,tab:grad:dft,tab:grad:matrix}). 
For this overview the notation of the gradients are changed from the explicit notation $\grad{z}$ to the shorter notation $\nabla_{z\conj}$. 
Also the alternative gradients $\nabla_{z} = \der{J}{z} = \left(\der{J}{z\conj}\right)\conj$ are shown.

\begin{table}[H]
	\rarray{1.3}
	\caption{Backward gradients of common scalar operations}\medskip
	\label{tab:grad:scalar}
	\renewcommand{\grad}[1]{\nabla_{#1\conj}}
	\newcommand{\gradz}[1]{\nabla_{#1}}
	\newcommand{\nSpalten}[2]{\multicolumn{#1}{ c }{#2}}
	\newcommand{\nZeilen}[2]{\multirow{#1}{*}{#2}}
	\newcommand{\linebreakcell}[2][l]{%
		\begin{tabular}[#1]{@{}#1@{}}#2\end{tabular}} 
	\robustify\bfseries
	\centering
	\begin{tabularx}{\textwidth}{llXX}
		\toprule
		 & Forward & Backward $\grad{\square}$ & Backward $\gradz{\square}$ \\
		\midrule
		Identity \eqref{eq:gradIdentity} & $s=z$ 
			& $\grad{z} = \gradIdentity$
			& $\gradz{z} = \gradz{s}$ 
			\\
		\midrule
		Conjugate \eqref{eq:gradConj} & $s=z\conj$ 
			& $\grad{z} = \gradConj$ 
			& $\gradz{z} = \gradz{s}\conj$ 
			\\
		\midrule
		Negation & $s=-z$ 
			& $\grad{z} = -\grad{s}$ 
			& $\gradz{z} = -\gradz{s}$ 
			\\
		\midrule
		\linebreakcell{Addition\\ \eqref{eq:gradAdd}} & $s=z_1+z_2$ 
			& \linebreakcell{$\grad{z_1} = \gradAdditionOne$\\$\grad{z_2} = \gradAdditionTwo $} 
			& \linebreakcell{$\gradz{z_1} = \gradz{s}$\\$\gradz{z_2} = \gradz{s}$} 
			\\
		\midrule
		\linebreakcell{Multiplication\\ \eqref{eq:gradMul}} & $s=z_1z_2$ 
			& \linebreakcell{$\grad{z_1} = \gradMultiplicationOne$\\$\grad{z_2} = \gradMultiplicationTwo$}
			& \linebreakcell{$\gradz{z_1} = \gradz{s}z_2$\\$\gradz{z_2} = \gradz{s}z_1$}  
			\\
		\midrule
		Exponentiation \eqref{eq:gradPow} & $s=z^n$ 
			& $\grad{z} = \gradPow$ 
			& $\gradz{z} = \gradz{s}n\left(z\right)^{n-1}$ 
			\\
		\midrule
		\linebreakcell{Division\\ \eqref{eq:gradDiv}} & $s=z_1/z_2$ 
			& \linebreakcell{$\grad{z_1} = \gradDivisionOne$\\$\grad{z_2} = \gradDivisionTwo$} 
			& \linebreakcell{$\gradz{z_1} = \gradz{s}/z_2$\\$\gradz{z_2} = -\gradz{s}\frac{z_1}{z_2^2}$}  
			\\
		\midrule
		Absolute \eqref{eq:gradAbs} & $s=\left|z\right|$ 
			& $\grad{z} = \gradAbsolute$ 
			& $\gradz{z} = \gradz{s}\e^{-\j\varphi}$ 
			\\

		\midrule
		Phase Factor & \nZeilen{2}{$s=\e^{\j\varphi}$} 
			& \multicolumn{2}{l}{$\grad{z} = \gradExpAngle$ } \\
			\eqref{eq:gradExpAngle} 
			& & \multicolumn{2}{r}{$\gradz{z} = \frac{\gradz{s}}{\left|z\right|} -\Re\left\{\gradz{s}\frac{1}{z\conj}\right\} \e^{-\j\varphi}$ } \\
		\midrule
		
		Real part \eqref{eq:gradReal} & $s=\Re\left\{z\right\}$ 
			& $\grad{z} = \gradReal$ 
			& $\gradz{z} = \gradz{s}$ 
			\\
		\midrule
		Imag. part \eqref{eq:gradImag} & $s=\Im\left\{z\right\}$ 
			& $\grad{z} = \gradImag$ 
			& $\gradz{z} = -\j\gradz{s}$ 
			\\
		\bottomrule
	\end{tabularx}
\end{table}

\begin{table}[H]
	\rarray{1.3}
	\caption{Backward gradients of discrete Fourier transformations}\medskip
	\label{tab:grad:dft}
	\renewcommand{\grad}[1]{\nabla_{#1\conj}}
	\newcommand{\gradz}[1]{\nabla_{#1}\hphantom{\conj}}
	\renewcommand{\s}{\vec{s}}
	\newcommand{\z}{\vec{z}}
	\renewcommand{\S}{\vec{S}}
	\newcommand{\Z}{\vec{Z}}
	\newcommand{\nZeilen}[2]{\multirow{#1}{*}{#2}}
	\newcommand{\linebreakcell}[2][l]{%
		\begin{tabular}[#1]{@{}#1@{}}#2\end{tabular}} 
	\robustify\bfseries
	\centering
	\begin{tabular}{lll}
		\toprule
		& \nZeilen{2}{Forward} & Backward $\grad{\square}$  \\
		\cline{3-3}	
		& & Backward $\gradz{\square}$ \\
		
		\midrule
		\linebreakcell{DFT} & \nZeilen{3}{$s_f = \dft{z_n}$ }
		& $\grad{z} = N \cdot \idft{\grad{s_f}}$ 
		\\
		\cline{3-3}
		\eqref{eq:grad:dft} & 
		& $\gradz{z} = \dft[f\rightarrow n]{\gradz{s_f}}$ 
		\\
		\midrule
		\linebreakcell{IDFT} & \nZeilen{3}{$s_f = \idft{z_n}$ }
		& $\grad{z} = \frac{1}{N} \cdot \dft{\grad{s_n}}$ 
		\\
		\cline{3-3}
		\eqref{eq:grad:idft} & 
		& $\gradz{z} = \idft[n\rightarrow f]{\gradz{s_n}}$ 
		\\
		\midrule
		\linebreakcell{RDFT} 
		\linebreakcell{\hphantom{A}\\\hphantom{A}\\\hphantom{A}\\\hphantom{A}} & \nZeilen{5}{$s_f = \rdft{z_n}$ }
		& {$\!\begin{aligned}
			\widetilde{\grad{s_f}}	&= \begin{cases}
			\Re\left\{\grad{s_f}\right\}, &\text{for } f\in\left\{0, \frac{N}{2}\right\} \\
			\frac{1}{2} \grad{s_f}, &\text{for } f\in\left[1, \dots, \frac{N}{2}-1\right]
			\end{cases} \\
			\grad{z_n} &= N \cdot \idft{\widetilde{\grad{s_f}}}
		\end{aligned}$}
		\\
		\cline{3-3}
		\linebreakcell[c]{\eqref{eq:grad:rdft1}\\and\\\eqref{eq:grad:rdft2}}
		\linebreakcell{\hphantom{A}\\\hphantom{A}\\\hphantom{A}\\\hphantom{A}} & 
		&  {$\!\begin{aligned}
			\widetilde{\gradz{s_f}}	&= \begin{cases}
			\Re\left\{\gradz{s_f}\right\}, &\text{for } f\in\left\{0, \frac{N}{2}\right\} \\
			\frac{1}{2} \gradz{s_f}, &\text{for } f\in\left[1, \dots, \frac{N}{2}-1\right]
			\end{cases} \\
			\gradz{z_n} &= \dft[f\rightarrow n]{\widetilde{\gradz{s_f}}}
			\end{aligned}$}
		\\
		\midrule
		\linebreakcell{IRDFT} 
		\linebreakcell{\hphantom{A}\\\hphantom{A}\\\hphantom{A}\\\hphantom{A}} & \nZeilen{5}{$s_f = \irdft{z_n}$ }
		& {$\!\begin{aligned}
		\widetilde{\grad{z_f}} &= \frac{1}{N} \cdot \dft{\grad{s_n}} \\
		\grad{z_f} &= \begin{cases}
		\widetilde{\grad{z_f}}, &\text{for } f \in\left\{0, \frac{N}{2}\right\} \\[0.2em]
		2\widetilde{\grad{z_f}}, &\text{for } f \in\left[1, \dots, \frac{N}{2}-1\right]
		\end{cases}
		\end{aligned}$}
		\\
		\cline{3-3}
		\linebreakcell[c]{\eqref{eq:grad:irdft1}\\and\\\eqref{eq:grad:irdft2}} \linebreakcell{\hphantom{A}\\\hphantom{A}\\\hphantom{A}\\\hphantom{A}} & 
		& {$\!\begin{aligned}
			\widetilde{\gradz{z_f}} &= \idft[n\rightarrow f]{\gradz{s_n}} \\
			\gradz{z_f} &=  \begin{cases}
			\widetilde{\gradz{z_f}}, &\text{for } f \in\left\{0, \frac{N}{2}\right\} \\[0.2em]
			2\widetilde{\gradz{z_f}}, &\text{for } f \in\left[1, \dots, \frac{N}{2}-1\right]
			\end{cases}
			\end{aligned}$}
		\\
		\bottomrule
	\end{tabular}
\end{table}

\begin{table}[H]
	\rarray{1.3}
	\caption{Backward gradients of common matrix operations}\medskip
	\label{tab:grad:matrix}
	\renewcommand{\grad}[1]{\nabla_{#1\conj}}
	\newcommand{\gradz}[1]{\nabla_{#1}\hphantom{\conj}}
	\renewcommand{\s}{\vec{s}}
	\newcommand{\z}{\vec{z}}
	\renewcommand{\S}{\vec{S}}
	\newcommand{\Z}{\vec{Z}}
	\newcommand{\nZeilen}[2]{\multirow{#1}{*}{#2}}
	\newcommand{\linebreakcell}[2][l]{%
		\begin{tabular}[#1]{@{}#1@{}}#2\end{tabular}} 
	\robustify\bfseries
	\centering
	\begin{tabular}{lll}
		\toprule
		& \nZeilen{2}{Forward} & Backward $\grad{\square}$  \\
		\cline{3-3}	
		& & Backward $\gradz{\square}$ \\
		
		\midrule
		\linebreakcell{Vector\\Normalization} \linebreakcell{\hphantom{A}\\\hphantom{A}} & \nZeilen{3}{$\s=\frac{\z}{\sqrt{\z\hermite\z}}$}
			& $\grad{\z} = \frac{\grad{\s}-\frac{\z}{\z\hermite\z}\Re\left\{\z\hermite\grad{\s}\right\}}{\sqrt{\z\hermite\z}}$ 
			\\
		\cline{3-3}
		\eqref{eq:grad:vecNormalisation} \linebreakcell{\hphantom{A}\\\hphantom{A}} & 
			& $\gradz{\z} = \frac{\gradz{\s}-\frac{\z\conj}{\z\hermite\z}\Re\left\{\z\transposed\gradz{\s}\right\}}{\sqrt{\z\hermite\z}}$
			\\
		
		\midrule
		\linebreakcell{Matrix\\Multiplication} & \nZeilen{3}{$\S=\Z_1\Z_2$} 
			& \linebreakcell{$\grad{\Z_1} = \grad{\S}\Z_2\hermite$ \\ $\grad{\Z_2} = \Z_1\hermite\grad{\S}$} 
			\\
		\cline{3-3}	
		\eqref{eq:grad:matMul1} and \eqref{eq:grad:matMul2} & 
			& \linebreakcell{$\gradz{\Z_1} = \gradz{\S}\Z_2\transposed$ \\ $\gradz{\Z_2} = \Z_1\transposed\gradz{\S}$} 
			\\
		
		\midrule
		\linebreakcell{Matrix Inversion} & \nZeilen{2}{$\S=\Z^{-1}$} 
			& $\grad{\Z} = -\S\hermite\grad{\S}\S\hermite$ 
			\\
		\cline{3-3}	
		\eqref{eq:grad:matInv} & 
			& $\gradz{\Z} = -\S\transposed\gradz{\S}\S\transposed$ 
			\\
		
		\midrule
		\linebreakcell{Matrix Inverse\\product 1} & \nZeilen{3}{$\S=\Z_1^{-1}\Z_{2}$} 
			& \linebreakcell{$\grad{\Z_2} = \Z_1\hermiteinv\grad{\S}$ \\ $\grad{\Z_1} = -\grad{\Z_2}\S\hermite$} 
			\\
		\cline{3-3}	
		\eqref{eq:grad:invMatMul11} and \eqref{eq:grad:invMatMul12} & 
			& \linebreakcell{$\gradz{\Z_2} = \Z_1\transposedinv\gradz{\S}$ \\ $\gradz{\Z_1} = -\gradz{\Z_2}\S\transposed$} 
			\\
		\midrule
		\linebreakcell{Matrix Inverse\\product 2} & \nZeilen{3}{$\S=\Z_1\Z_2^{-1}$} 
			& \linebreakcell{$\grad{\Z_1} = \grad{\S}\Z_2\hermiteinv$ \\ $\grad{\Z_2} = -\S\hermite\grad{\Z_1}$} 
			\\
		\cline{3-3}	
		\eqref{eq:grad:invMatMul21} and \eqref{eq:grad:invMatMul22} & 
			& \linebreakcell{$\gradz{\Z_1} = \gradz{\S}\Z_2\transposedinv$ \\ $\gradz{\Z_2} = -\S\transposed\gradz{\Z_1}$} 
			\\
		
		\midrule
		\linebreakcell{Cholesky\\decomposition} \linebreakcell{\hphantom{A}\\\hphantom{A}} & \nZeilen{3}{$\S=\vec{g}\left(\Z\right)$} 
			& \linebreakcell{$\grad{\tilde{\Z}} = \S\hermiteinv\left(\left(\S\hermite\grad{\S}\right)\hadamard \left(\lltriangle[very thick] - \frac{1}{2}\vec{I}\right) \right)\S\inv$\\
			$\grad{\Z} = \frac{1}{2}\left(\grad{\tilde{\Z}} + \left(\grad{\tilde{\Z}}\right)\hermite\right)$}
			\\
		\cline{3-3}	
		\eqref{eq:chol:grad} and \eqref{eq:chol:gradStable} \linebreakcell{\hphantom{A}\\\hphantom{A}}  & 
			& \linebreakcell{$\gradz{\tilde{\Z}} = \S\transposedinv\left(\left(\S\transposed\gradz{\S}\right)\hadamard \left(\lltriangle[very thick] - \frac{1}{2}\vec{I}\right) \right)\left(\S\inv\right)\conj$\\
				$\gradz{\Z} = \frac{1}{2}\left(\gradz{\tilde{\Z}} + \left(\gradz{\tilde{\Z}}\right)\hermite\right)$}
			\\
		
		\midrule
		\linebreakcell{Eigenvalue\\decomposition} & \nZeilen{4}{$\vecG{\Lambda}, \vec{\S}= \mathrm{eig}\left(\Z\right)$} & 
			
			\multicolumn{1}{l}{$\!\begin{aligned}
				\grad{\Z} &= \S\hermiteinv \left(\grad{\vecG{\Lambda}} + \vec{F}\conj \circ \left(\S\hermite \grad{\S} \right)\right) \S\hermite \\
				&\hphantom{=} -\S\hermiteinv \left(\vec{F}\conj \circ \S\hermite\S\left(\Re\{\S\hermite\grad{\S}\}\circ\vec{I}\right)\right) \S\hermite \\
				F_{ij} &= \left(\lambda_j - \lambda_i\right)\inv \quad \forall_{i \ne j} \quad\text{and}\quad F_{ii} = 0 
			\end{aligned}$}
			\\
		\cline{3-3}	
		\eqref{eq:eig:final}\linebreakcell{\hphantom{A}\\\hphantom{A}} & &
			\multicolumn{1}{l}{$\!\begin{aligned}
				\gradz{\Z} &= \S\transposedinv \left(\gradz{\vecG{\Lambda}} + \vec{F}\hphantom{\conj} \circ \left(\S\transposed\gradz{\S} \right)\right) \S\transposed \\
				&\hphantom{=} -\S\transposedinv \left(\vec{F}\hphantom{\conj} \circ \S\transposed\S\conj\left(\Re\{\S\transposed\gradz{\S}\}\circ\vec{I}\right)\right) \S\transposed
				\end{aligned}$}
		 \\
		\bottomrule
	\end{tabular}
\end{table}


\chapter{Evaluation} 
\label{cha:evaluation}

To evaluate the approach described in \cref{cha:status} two different performance measures are used.
The quality of the beamformed signal is evaluated with the \gls{pesq} measure.
Also the \gls{snr} is calculated.
Both measures are described in the following sections.
Afterwards the system environment and hyperparameters as well as the simulation procedure are denoted.
In the end an outlook to an end-to-end system evaluated with the \gls{wer} is given.

\section{Perceptual Evaluation of Speech Quality (PESQ)}

In the ITU-T-Recommendation P.862 (International Telecommunication Union) and \cite{Rix2001PESQ} the \gls{pesq} procedure for wideband signals is described.
It compares a transmitted with a received signal via an psychoacoustic model.
The measurement is done in MOS-LQO (Mean Opinion Score - Listening Quality Objective).
The MOS-Scale is defined as follows:
$5=\text{Excellent}$, $4=\text{Good}$, $3=\text{Fair}$, $2=\text{Poor}$, $1=\text{Bad}$.

\begin{figure}[H]
	\centering
	\begin{tikzpicture}
		\begin{axis}[
			hide axis,
			scale only axis,
			height=0pt,
			width=0pt,
			colormap/jet,
			colorbar horizontal,
			point meta min=1,
			point meta max=5,
			colorbar style={
				width=10cm,
				xtick={5,4,3,2,1},
			}]
			\addplot [draw=none] coordinates {(0,0)};
			\end{axis}
			\begin{axis}[
			hide axis,
			scale only axis,
			height=0pt,
			width=0pt,
			colormap/jet,
			colorbar horizontal,
			point meta min=1,
			point meta max=5,
			colorbar style={
				axis on top,
				width=10cm,
				xtick={5,4,3,2,1},
				xticklabels={{Excellent},{Good},{Fair},{Poor},{Bad}},
				xticklabel pos=right
			}]
			\addplot [draw=none] coordinates {(0,0)};
		\end{axis}
	\end{tikzpicture} \medskip
	\caption{MOS-Scale}
\end{figure}
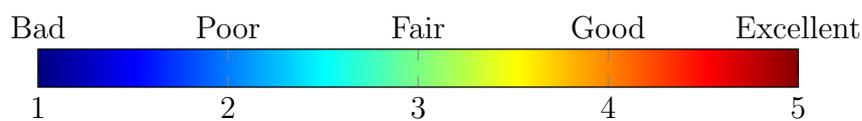

\section{Signal-to-noise ratio (SNR)}

The \gls{snr} Criterion describes the ratio between the useful signal power and the noise power. The input \gls{snr} 
\begin{align}
	\mathrm{SNR}_\mathrm{in} &= 10 \log_{10} \frac{ \sumT \sumF \left|\left|\vec{X}_{f,t}\right|\right|^2 }{\sumT \sumF \left|\left|\vec{N}_{f,t}\right|\right|^2}	
\end{align}
is the ratio between the clean signal and the noise and the output \gls{snr} in a beamformer scenario
\begin{align}
	\mathrm{SNR}_\mathrm{out} &= 10 \log_{10} \frac{ \sumT \sumF \left|\vec{w}_f\hermite\vec{X}_{f,t}\right|^2 }{\sumT \sumF \left|\vec{w}_f\hermite\vec{N}_{f,t}\right|^2}	
\end{align}
the ratio of the beamformed clean signals and the beamformed noise.

\section{Environment and hyperparameter}

\begin{itemize}
	\item Database: CHiME-3 \cite{Barker2015Chime}
	\begin{itemize}
		\item Sets:
		\begin{itemize}
			\item Train: \textit{tr05\textunderscore simu} with \num{7138} simulated utterances
			\item cross validation: \textit{dt05\textunderscore simu} with \num{1640} simulated utterances
			\item evaluation: \textit{dt05\textunderscore simu} with \num{1640} simulated utterances
		\end{itemize}
		\item Microphones:
		\begin{itemize}
			\item Count: \num{6}
			\item Positions: fixed
			\item Speaker positions: low variation 
		\end{itemize}
		\item Locations:
		\begin{itemize}
			\item bus
			\item cafe
			\item pedestrian area
			\item  street junction
		\end{itemize}
		\item Text prompts: 
		\begin{itemize}
			\item Language: English
			\item Origin: Wall Street Journal (WSJ0) corpus
		\end{itemize}
	\end{itemize}
	\item Sampling frequency: \SI{16000}{Hz}
	\item \gls{stft} parameters
	\begin{itemize}
		\item FFT length: \num{1024} ($=\SI{64}{ms}$
		\item FFT shift: \num{256} ($=\SI{16}{ms}$, $\frac{1}{4}$ of the length)
		\item window: blackman
	\end{itemize}

\end{itemize}

For the evaluation of the \gls{NN} we used the same Dataset as for cross validation.
The official evaluation dataset of CHiME-3 can not be used, because it does not contain the clean signal, which is required for the performance measurements.

\section{Simulation}
To evaluate the complex gradients, they are tested against \cite{Heymann2016}, 
which uses the \gls{NN} to obtain masks for a subsequent beamforming step (see \cref{tab:config}).
In \cite{Heymann2016} the objective function is the binary cross entropy of the network output with a heuristic binary mask. 
The objective function in this report is replaced with an \gls{snr} objective, like described in \eqref{eq:snr:loss}.

\cref{tab:scores} summarizes the simulation results. 
The first row shows the input \gls{pesq} and \gls{snr}. 
The first column indicates the origin of the masks, where Oracle \cite{Heymann2016} is the heuristic binary mask,
\cite{Heymann2016} the network trained with the binary cross entropy objective and NN-GEV and NN-MVDR use mask estimation networks trained with an \gls{snr} objective,
where the first uses a \gls{gev} beamformer and the second an \gls{mvdr} beamformer. 
Note that the beamformer used during \gls{NN} training (denoted "mask source" in \cref{tab:scores}) can be different from the evaluation beamformer. 
The optional post filter \gls{ban} is described in \cite{warsitz2007ban} and $||\cdot||_2$ is a normalization to unit length.

The proposed alternative objective with a GEV beamformer outperforms the reference \cite{Heymann2016} both in \gls{pesq} and \gls{snr}. 
Unexpectedly the \gls{gev} beamformer (also known as max \gls{snr} beamformer) nearly reaches the oracle perfomance in \gls{pesq}, but not the oracle \gls{snr}.
The \gls{mvdr} beamformer as objective is not able to beat the reference with the \gls{mvdr} evaluation beamformer.
With an \gls{gev} evaluation the results get better, but they do not exceed the NN-GEV.

\begin{table}
	\rarray{1.3}
	\caption{Configuration of the neural network}\medskip
	\label{tab:config}
	\robustify\bfseries
	\centering
	\begin{tabular}{ccccc}
		\toprule
		& Units & Type & Non-Linearity & $p_\mathrm{dropout}$ \\
		\midrule
		L1 & 256 & BLSTM & Tanh & $0.5$ \\
		L2 & 513 & FF & ReLU & $0.5$ \\
		L3 & 513 & FF & ReLU & $0.5$ \\
		L4 & 1026 & FF & Sigmoid & $0.0$ \\
		\bottomrule
	\end{tabular}
\end{table}
In \cref{fig:masks} some masks at the output of the \gls{NN} are shown exemplary, 
where \cref{fig:masks:bce,fig:masks:gev} are from \cite{Heymann2016} and NN-GEV, respectively.
The evaluation beamformer for the scores is a GEV with a BAN postfilter in both cases.
Counter-intuitively the right mask yields slightly better scores, 
although no expected symmetries as in the left one are visible.
Also the mask tends to suddenly change its values between the frequencies.
{
\newcommand{\NUM}[1]{\num[round-mode=places,round-precision=1]{#1}}
\newcommand{\snr}{SNR:}
\newcommand{\pesq}{PESQ:}
\newcommand{\uttID}{f04_051c0112_str}
\begin{figure}
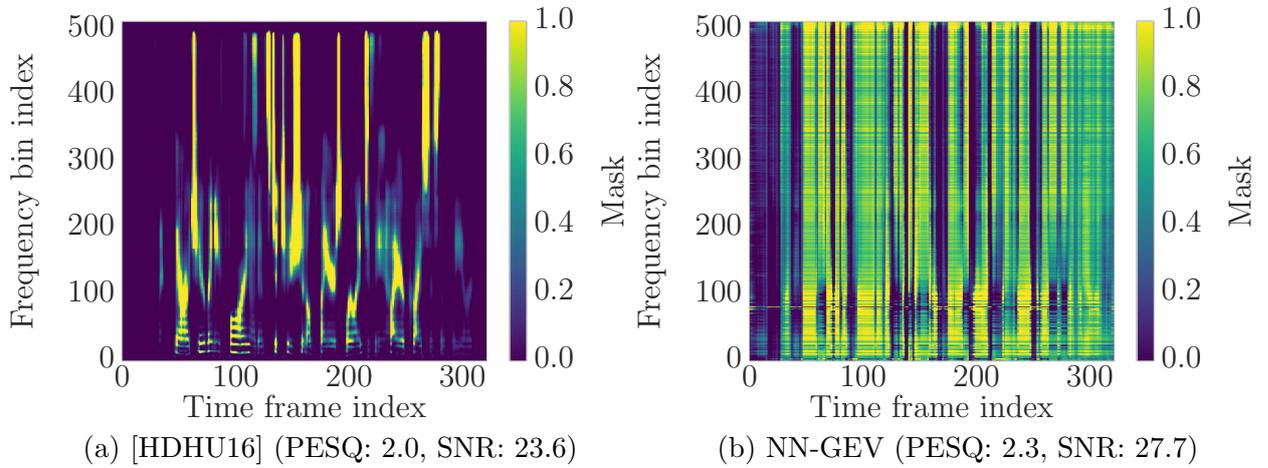

	\subcaptionbox{
		\cite{Heymann2016} (\gls{pesq}: \NUM{2.013833}, \gls{snr}: \NUM{23.558835})\label{fig:masks:bce}
	}{
		\def\svgwidth{0.5\textwidth}
		\import{images/bce/}{\uttID_X_mask.eps_tex}
	}
	\subcaptionbox{
		NN-GEV (\gls{pesq}: \NUM{2.304666}, \gls{snr}: \NUM{27.679976})\label{fig:masks:gev}
	}{
		\def\svgwidth{0.5\textwidth}
		\import{images/gev/}{\uttID_X_mask.eps_tex}
	}
	\caption[Sample speaker output masks for reference and NN-GEV]{Sample output masks for reference and NN-GEV\footnotemark}
	\label{fig:masks}
\end{figure}
\footnotetext{The CHiME-3 utterance ID is \texttt{\detokenize{f04_051c0112_str}}.}
}

\begin{table}[!h]
	\rarray{1.3}
	\caption[Overview of the results for different system configurations.]{Overview of the results for different system configurations. The first row shows the scores of the input noisy signal.} \medskip
	\label{tab:scores}
	\robustify\bfseries
	\centering
	\newcommand{\myRot}[1]{\rotatebox[origin=c]{90}{#1}}
	\newcommand{\ban}{BAN}
	\newcommand{\BCE}{\cite{Heymann2016}}
	\newcommand{\comment}[1]{}
	\newcommand{\NUM}[1]{\num[round-mode=places,round-precision=2]{#1}}
	\newcommand{\nRwo}[2]{\multirow{#1}{*}{#2}}
	\newcommand{\linebreakcell}[2][c]{%
		\begin{tabular}[#1]{@{}c@{}}#2\end{tabular}} 
	\begin{tabular}{ccccc}
		\toprule
		{\linebreakcell{mask\\source}} & {\linebreakcell{evaluation\\beamformer}} & {\linebreakcell{post\\filter}} & PESQ & SNR \\ 
		\midrule
		-					& -					& -			& \NUM{1.211}		& \NUM{3.850991}		\\ 
		\midrule
		\nRwo{2}{Oracle \cite{Heymann2016}}& \nRwo{2}{GEV}	& -	& \NUM{1.725129}& \NUM{15.349016}		\\ 
		& 	& \ban		& \NUM{1.764146}	& \NUM{18.016347}		\\ 
		\midrule
		\nRwo{2}{\BCE}	& \nRwo{2}{GEV}			& -			& \NUM{1.640623}	& \NUM{14.078099}	\\ 
		& 					& \ban		& \NUM{1.699941}	& \NUM{16.029049}	\\ 
		\midrule
		\nRwo{2}{NN-GEV}	& \nRwo{2}{GEV}		& -			& \NUM{1.713171}	& \NUM{14.588011}	\\ 
		& 					& \ban		& \NUM{1.746576}	& \NUM{16.045985}	\\ 
		\midrule
		\nRwo{4}{NN-MVDR}	& \nRwo{2}{MVDR}	& -			& \NUM{1.434078}	& \NUM{10.575671}	\\ 
		&			 	& $||\cdot||_2$	& \NUM{1.596658}	& \NUM{13.843453}	\\ 
		\cline{2-5}	
		& \nRwo{2}{GEV}		& -			& \NUM{1.598684}	& \NUM{14.211392}	\\ 
		& 					& \ban		& \NUM{1.502961}	& \NUM{14.511593}	\\ 
		\bottomrule
	\end{tabular}
\end{table}

\FloatBarrier

\section{Outlook: ASR end-to-end}
\label{sec:eval:endToEnd}
Based on this work in \cite{heymann2017EndToEnd} further simulations with a backpropagation through beamformers are done.
Here only a small summary is presented.
There the larger system (see \cref{fig:blockschaltbild} purple part) is described with two \glspl{NN} in an \gls{asr} system. 
The two \glspl{NN} are one beamforming \gls{NN} (BF), which is the same as described here, and an acoustic model \gls{NN} (AM).
\cref{tab:all_systems_average} shows the \gls{wer} from \cite{heymann2017EndToEnd}.

The training is done with different configurations:
\begin{itemize}
	\item fixed: Each network is trained separately, without knowing the other \gls{NN}. This means that the AM is trained on each channel of the noisy signal.
	\item scratch: Networks are randomly initialized and only trained in an end-to-end configuration, where the gradient of the AM is passed to the mask estimation \gls{NN}.
	\item finetune: The networks are first trained separately an then finetuned in an end-to-end configuration.
\end{itemize}

\begin{table}[hbt]
	\centering
	\caption[Average \gls{wer} (\%) for the described systems in \cite{heymann2017EndToEnd}.]{Average \gls{wer} (\%) for the described systems in \cite{heymann2017EndToEnd}. Copied with permission from \cite{heymann2017EndToEnd}.}
	\label{tab:all_systems_average}
	\rarray{1.3}
	\definecolor{mygray}{RGB}{0,0,0}
	\setlength{\tabcolsep}{2pt}
	\begin{tabular}{C{2cm}C{2cm}C{1cm}C{1cm}C{1cm}C{1cm}}
		\toprule
		\multicolumn{2}{c}{Training} & \multicolumn{2}{c}{Dev} & \multicolumn{2}{c}{Test}\\
		\midrule
		BF & AM & real & simu & real & simu\\
		\midrule
		\textcolor{mygray}{fixed} & \textcolor{mygray}{fixed} & \textcolor{mygray}{4.26} & \textcolor{mygray}{4.29} & \textcolor{mygray}{5.85} & \textcolor{mygray}{4.59} \\
		scratch & scratch & 5.51 & 5.19 & 8.76 & 5.61 \\
		scratch & finetune & 4.14 & 4.09 & 5.86 & 4.06 \\
		fixed & finetune & 4.09 & 3.96 & 5.56 & {3.9} \\
		finetune & finetune & {3.77} & {3.89} & {5.42} & 3.95 \\
		\bottomrule
	\end{tabular}
\end{table}

While pure end-to-end training was not able to increase the performance, finetuning the backprobagation through the beamformer outperforms the reference system.


\chapter{Summary} 
\label{cha:summary}

In this report the derivation of \cv gradients used in \glspl{NN} is described.
Many well-known \rv gradients are extended to their \cv counterparts.
This enables the use of \cv operations in the calculation of the objective function of a \gls{NN}.
Here the optimization of a mask estimation \gls{NN} \wrt the beamformer output \gls{snr} is shown exemplarily.
The changed objective function leads to a gain in \gls{pesq} and \gls{snr} in comparison to a objective function based on a heuristic reference mask.
Also two different statistically optimal beamformer operations are compared.

As a key theoretical result the derivation of the gradient of the eigenvalue problem with \cv eigenvalues is derived,
which is an extension of the known gradient for hermitian matrices.
Also the \cv gradients of the Cholesky decomposition, matrix product and inverse, \gls{dft} and vector normalization are derived.
So this report can be treated as a collection of \cv gradients, which can be used in the calculation of a \gls{NN} objective function, too.

Using these results the complete \gls{asr} system with a mask estimation \gls{NN} and a multichannel beamformer as a speech enhancement preprocessing can be optimized together \wrt to a common objective.
The results of simulations described in \cref{sec:eval:endToEnd} show that pretrained networks fine-tuned with the backpropagation trough the beamformer leads to better \glspl{wer}.

    \appendix

\chapter{Appendix} 
\label{cha:appendix}

\section{Proof of the matrix product rule}
\label{subsec:matrixproduct}
\begin{align}
\derz{AB} = \derz{A}B + A\derz{B} \qquad \zeta \in {z, z\conj}
\end{align}
\begin{proof}
	\newcommand*{\derx}[1]{\der{#1}{x}}
	\newcommand*{\dery}[1]{\der{#1}{y}}
	\renewcommand*{\derz}[1]{\der{#1}{z\conj}}
	With $C,D,G,H \in \Real$ and
	\begin{align}
	A = C + \j D, \quad B = G + \j H, \quad AB = CG - DH + \j \left(DG + CH \right)
	\end{align}
	\begin{align}\begin{aligned}
	\der{AB}{\zeta} =& \frac{1}{2} \left(\derx{AB} \pm \j \dery{AB}\right) = \begin{cases}
	\frac{1}{2} \left(\derx{AB} + \j \dery{AB}\right) &\text{for } \zeta = z\conj\\
	\frac{1}{2} \left(\derx{AB} - \j \dery{AB}\right) &\text{for } \zeta = z
	\end{cases}  \\
	=& \frac{1}{2} \left(\derx{}\left(CG-DH\right) \mp \dery{}\left(DG+CH\right)+\j\left(\derx{}\left(DG+CH\right)\pm\dery{}\left(CG-DH\right)\right)\right)\\
	=& \frac{1}{2}\left(\derx{C}G + C\derx{G} - \derx{D}H - D\derx{H} \mp \dery{D}G \mp D\dery{G} \mp \dery{C}H \mp c\dery{H} \right) \\
	&\pm \frac{1}{2}\j \left(\derx{D}G + D\derx{G} + \derx{C}H+C\derx{H}\pm\dery{C}G\pm C\dery{G}\mp\dery{D}H\mp D\dery{H} \right) \\
	=& \frac{1}{2}\left(\derx{C}\underbrace{(G+\j H)}_{B} - \derx{D}\underbrace{(H-\j G)}_{-\j B} \mp \dery{D}\underbrace{(G+\j H)}_{B} \mp \dery{C}\underbrace{(H-\j G)}_{-\j B}\right) \\
	&+ \frac{1}{2}\left(\underbrace{(C+\j D)}_A\derx{G} -\underbrace{(D-\j C)}_{-\j A}\derx{H} \mp \underbrace{(D-\j C)}_{-\j A} \dery{G} \mp\underbrace{(C+\j D)}_A\dery{H} \right) \\
	=& \frac{1}{2}\left(\derx{C}+\j\derx{D}\mp\dery{D}\pm\j\dery{C}\right)B + \frac{1}{2}A\left(\derx{G}+\j\derx{H}\mp\dery{H}\pm\j\dery{G}\right) \\
	=& \der{A}{\zeta}B + A \der{B}{\zeta}
	\end{aligned}\end{align}
\end{proof}

\section{Proof of complex-valued chain rule}
\label{appendix:proof:complex_valued_chain_rule}

\begin{proof}
	For the verification of \eqref{eq:chainRule} 
\begin{align}
	\der{J}{z\conj} &=\frac{1}{2} \left(\der{J}{\sigma} \der{\sigma}{x} + \der{J}{\omega} \der{\omega}{x}+ \j \der{J}{\sigma} \der{\sigma}{y} + \j \der{J}{\omega} \der{\omega}{y}\right)\nonumber\\
	&\stackrel{!}{=} \left(\der{J}{s\conj}\right)\conj \der{s}{z\conj} + \der{J}{s\conj} \left(\der{s}{z}\right)\conj, \nonumber
\end{align}
	it is easier to start with the right side of the equation, inserting \eqref{eq:Wirtinger} and simplify it. This leads to: 
\begin{align}\begin{aligned}
	&4\left(\der{J}{s\conj}\right)\conj\der{s}{z\conj} + 4\der{J}{s\conj}\left(\der{s}{z\conj}\right) \\
	&= \frac{4}{4} \left(\der{J}{\sigma}-\j\der{J}{\omega}\right) \left(\der{s}{x}+\j\der{s}{y}\right)
	+ \frac{4}{4} \left(\der{J}{\sigma}+\j\der{J}{\omega}\right) \left(\der{s\conj}{x}+\j\der{s\conj}{y}\right), \\
	&= \der{J}{\sigma} \der{s}{x}+\der{J}{\omega} \der{s}{y}+j \der{J}{\sigma} \der{s}{y}-\j \der{J}{\omega} \der{s}{x} 
	+ \der{J}{\sigma} \der{s\conj}{x} - \der{J}{\omega} \der{s\conj}{y} +\j \der{J}{\sigma} \der{s\conj}{y} +\j \der{J}{\omega} \der{s\conj}{x},  \\
	&= 
	\der{J}{\sigma} \left(\der{s}{x} + \der{s\conj}{x} \right)
	+\der{J}{\omega} \left(\der{s}{y} - \der{s\conj}{y} \right) 
	+\j \der{J}{\sigma} \left(\der{s}{y} + \der{s\conj}{y} \right)
	-\j \der{J}{\omega} \left(\der{s}{x} - \der{s\conj}{x} \right),\\
	&=
	\der{J}{\sigma} \der{}{x}\underbrace{\left(s + s\conj\right)}_{2\sigma}
	+\der{J}{\omega} \der{}{y} \underbrace{\left(s - s\conj \right)}_{\j 2\omega}
	+\j \der{J}{\sigma} \der{}{y} \underbrace{\left(s + s\conj \right)}_{2\sigma}
	-\j \der{J}{\omega} \der{}{x} \underbrace{\left(s - s\conj \right)}_{\j 2\omega},
	\\
	&= 2 \left(
	\der{J}{\sigma} \der{\sigma}{x}
	+\der{J}{\omega} \der{\omega}{x}
	+\j \der{J}{\sigma} \der{\sigma}{y}
	+\j \der{J}{\omega} \der{\omega}{y}
	\right).
\end{aligned}\end{align}
\end{proof}

\pagebreak
\section{Verification of a matrix identity for the calculation of the eigendecomposition gradient}
\label{appendix:eigen:CLambdaMinusLamdaC}
We want to prove:
\begin{align}
	\vec{C} \vecG{\Lambda} - \vecG{\Lambda} \vec{C} &= \vec{C} \hadamard \vec{E}, \\
\end{align}
where
\begin{align}
	\vecG{\Lambda} &= \left[\begin{array}{ccc}
			\lambda_1 & 0 & \dots\\
			0 & \lambda_2 & \\
			\vdots &  & \ddots
		\end{array}\right] \\
	\vec{C} &= \left[\begin{array}{ccc}
			c_{11} & c_{21} & \dots\\
			c_{21} & c_{22} & \\
			\vdots &  & \ddots
		\end{array}\right] \\
	\vec{E} &= \left[\begin{array}{ccc}
			\lambda_1 - \lambda_1 & \lambda_2 - \lambda_1 & \dots\\
			\lambda_1 - \lambda_2 & \lambda_1 - \lambda_1 & \\
			\vdots &  & \ddots \\
		\end{array}\right] 
		= \left[\begin{array}{ccc}
			E_{11} & E_{12} & \dots\\
			E_{21} & E_{22} & \\
			\vdots &  & \ddots \\
		\end{array}\right] \\
	E_{ij} &= \lambda_j - \lambda_i 
\end{align}

\begin{proof}
	\begin{align}
		\vec{C} \vecG{\Lambda} - \vecG{\Lambda} \vec{C} &= 
			\left[\begin{array}{ccc}
				c_{11} & c_{21} & \dots\\
				c_{21} & c_{22} & \\
				\vdots &  & \ddots
			\end{array}\right]
			\left[\begin{array}{ccc}
				\lambda_1 & 0 & \dots\\
				0 & \lambda_2 & \\
				\vdots &  & \ddots
			\end{array}\right]
			- \left[\begin{array}{ccc}
				\lambda_1 & 0 & \dots\\
				0 & \lambda_2 & \\
				\vdots &  & \ddots
			\end{array}\right]
			\left[\begin{array}{ccc}
				c_{11} & c_{21} & \dots\\
				c_{21} & c_{22} & \\
				\vdots &  & \ddots
			\end{array}\right]
		 \nonumber\\
		 &= 
			 \left[\begin{array}{ccc}
			 	\lambda_1 c_{11} & \lambda_2 c_{21} & \dots\\
			 	\lambda_1 c_{21} & \lambda_2 c_{22} & \\
			 	\vdots &  & \ddots
			 \end{array}\right]
			 - \left[\begin{array}{ccc}
			 	\lambda_1 c_{11} & \lambda_1 c_{21} & \dots\\
			 	\lambda_2 c_{21} & \lambda_2 c_{22} & \\
			 	\vdots &  & \ddots
			 \end{array}\right]
		 \nonumber\\
		 &= 
		\vec{C}\hadamard\left[\begin{array}{ccc}
		 	\lambda_1 - \lambda_1 & \lambda_2 - \lambda_1 & \dots\\
		 	\lambda_1 - \lambda_2 & \lambda_2 - \lambda_2 & \\
		 	\vdots &  & \ddots \\
		 \end{array}\right]  \nonumber\\
		 &= \vec{C}\hadamard\vec{E} \nonumber
	\end{align}
\end{proof}

\section{Gradient identity proof for the eigenvalue decomposition}
\label{appendix:proof:eig_identity}

The operation
\begin{align}
	\vec{\Lambda}, \vec{W} &= \mathrm{eig}(\vecG{\phi}), \\
	\tilde{\vecG{\phi}} &= \vec{W} \vec{\Lambda} \vec{W}\inv,
\end{align}
does not effect the gradient 
\begin{align}
	\grad{\vecG{\phi}} \stackrel{!}{=} \grad{\tilde{\vecG{\phi}}} .
\end{align}

\begin{proof}
	We start by introducing some auxiliary variables
	\begin{align}
		\tilde{\vecG{\phi}} &= \vec{W}_1 \vec{\Lambda} \vec{W}_2\inv = \vec{W}_1 \vec{M}. \\
	\end{align}
	Using the gradients for matrix multiplication \eqref{eq:grad:matMul1} and \eqref{eq:grad:matMul2})
	\begin{align}
		\grad{\vec{W}_1} &= \grad{\tilde{\vecG{\phi}}} \vec{M}\hermite 
		= \grad{\tilde{\vecG{\phi}}} \vec{W}\hermiteinv\vec{\Lambda}\hermite, \\
		\grad{\vec{M}} &= \vec{W}_1\hermite \grad{\tilde{\vecG{\phi}}}, \\
	\end{align}
	the gradient of the matrix inverse multiplication \eqref{eq:grad:invMatMul21} and \eqref{eq:grad:invMatMul22} under the constraint that $\grad{\vec{\Lambda}}$ is a diagonal matrix
	\begin{align}
		\grad{\vec{\Lambda}} &= \left(\grad{\vec{M}}\vec{W}_2\hermiteinv\right) \hadamard \vec{I} 
		= \left(\vec{W}\hermite \grad{\tilde{\vecG{\phi}}}\vec{W}\hermiteinv\right) \hadamard \vec{I}, \\
		\grad{\vec{W}_2} &= -\vec{M}\hermite\grad{\vec{M}}\vec{W}_2\hermiteinv
		= -\vec{W}\hermiteinv\vec{\Lambda}\hermite\vec{W}\hermite \grad{\tilde{\vecG{\phi}}}\vec{W}\hermiteinv,
	\end{align}
	the rule for reusing a variable (\cref{sec:reuse_variable})
	\begin{align}
		\grad{\vec{W}} &= \grad{\vec{W}_1} + \grad{\vec{W}_2} 
		= \grad{\tilde{\vecG{\phi}}} \vec{W}\hermiteinv\vec{\Lambda}\hermite 
		- \vec{W}\hermiteinv\vec{\Lambda}\hermite\vec{W}\hermite \grad{\tilde{\vecG{\phi}}}\vec{W}\hermiteinv,
	\end{align}
	the gradient for the eigenvalue decomposition \eqref{eq:eig:backwardgradient}
	\begin{align}
		\grad{\tilde{\vecG{\phi}}} &= \vec{W}\hermiteinv \left(\grad{\Lambda}+\vec{F}\conj\hadamard\left(\vec{W}\hermite\grad{\vec{W}}\right)\right) \vec{W}\hermite \\
		&= \vec{W}\hermiteinv \left(\grad{\Lambda}+\vec{F}\conj\hadamard\left(\vec{W}\hermite\grad{\tilde{\vecG{\phi}}} \vec{W}\hermiteinv\vec{\Lambda}\hermite 
		- \vec{W}\hermite\vec{W}\hermiteinv\vec{\Lambda}\hermite\vec{W}\hermite \grad{\tilde{\vecG{\phi}}}\vec{W}\hermiteinv\right)\right) \vec{W}\hermite \\
		&= \vec{W}\hermiteinv \left(\grad{\Lambda} +\vec{F}\conj\hadamard\left(\vec{W}\hermite\grad{\tilde{\vecG{\phi}}} \vec{W}\hermiteinv\vec{\Lambda}\hermite 
		- \vec{\Lambda}\hermite\vec{W}\hermite \grad{\tilde{\vecG{\phi}}}\vec{W}\hermiteinv\right)\right) \vec{W}\hermite
	\end{align}
	and considering that $\vec{C}\vec{\Lambda} - \vec{\Lambda}\vec{C} = \vec{E} \hadamard \vec{C}$ \eqref{eq:eig:tmp1} and $\vec{E}\hadamard\vec{F} = \vec{1} -\vec{I}$ \eqref{eq:eig:tmp3}
	\begin{align}
		\vec{C}\vec{\Lambda}\hermite - \vec{\Lambda}\hermite\vec{C} 
		&= \vec{C}\vec{\Lambda}\conj - \vec{\Lambda}\conj\vec{C} 
		= \vec{E}\conj \hadamard \vec{C}\\
		\grad{\tilde{\vecG{\phi}}}
		&= \vec{W}\hermiteinv \left(\grad{\Lambda} +\vec{F}\conj\hadamard \vec{E}\conj\hadamard\left(\vec{W}\hermite\grad{\tilde{\vecG{\phi}}} \vec{W}\hermiteinv \right)\right) \vec{W}\hermite \\
		&= \vec{W}\hermiteinv \left(\left(\vec{W}\hermite \grad{\tilde{\vecG{\phi}}}\vec{W}\hermiteinv\right) \hadamard \vec{I} +\left(\vec{1}-\vec{I}\right)\hadamard\left(\vec{W}\hermite\grad{\tilde{\vecG{\phi}}} \vec{W}\hermiteinv \right)\right) \vec{W}\hermite \\
		&= \vec{W}\hermiteinv \left(\left(\vec{W}\hermite \grad{\tilde{\vecG{\phi}}}\vec{W}\hermiteinv\right) \hadamard \left(\vec{I} + \vec{1}-\vec{I}\right)  \right) \vec{W}\hermite \\
		&= \grad{\tilde{\vecG{\phi}}}
	\end{align}
	verifies that the gradient is not effected by this operation, under the assumption that the gradient extension of eigenvalue decomposition is zero.
	The assumption holds since 
	\begin{align}
		\Re\left\{\vec{W}\hermite\grad{\vec{W}}\right\} 
		&= \Re\left\{\vec{W}\hermite\grad{\tilde{\vecG{\phi}}} \vec{W}\hermiteinv\vec{\Lambda}\hermite 
		- \vec{W}\hermite\vec{W}\hermiteinv\vec{\Lambda}\hermite\vec{W}\hermite \grad{\tilde{\vecG{\phi}}}\vec{W}\hermiteinv\right\},
	\end{align}
	with $\vec{I} = \vec{W}\hermite\vec{W}\hermiteinv$ and $\vec{C}\vec{\Lambda}\hermite - \vec{\Lambda}\hermite\vec{C} = \vec{E}\conj \hadamard \vec{C}$
	\begin{align}
		\Re\left\{\vec{W}\hermite\grad{\vec{W}}\right\} 
		&= \Re\left\{\vec{W}\hermite\grad{\tilde{\vecG{\phi}}} \vec{W}\hermiteinv\vec{\Lambda}\hermite 
		- \vec{\Lambda}\hermite\vec{W}\hermite \grad{\tilde{\vecG{\phi}}}\vec{W}\hermiteinv\right\} \\
		&= \Re\left\{\left(\vec{W}\hermite\grad{\tilde{\vecG{\phi}}} \vec{W}\hermiteinv\right) \hadamard \vec{E}\conj \right\} 
	\end{align}
	is a off diagonal matrix (the elements on the diagonal are zero) and appears in the extended gradient $\Re\left\{\vec{W}\hermite\grad{\vec{W}}\right\} \hadamard \vec{I}$. Therefore the extended gradient is for this example not necessary.
\end{proof}

\section{GEV beamformer via spatial whitening}
\label{appendix:GEV_spatial_whitening}

Noise \gls{PSD} Matrix:
\begin{align}
{\vecG{\Phi}_{NN}}_f &= \frac{\sumT\sumD M^{(\nu)}_{f, t, d}\Y_{f, t} \cdot \Y_{f, t}\hermite}{\sumT \sumD M^{(\nu)}_{f, t, d}}.
\end{align}	

Spatial Whitening:
\begin{align}
	\Ytilde_{f, t} &= \phinn_{f}^{-\frac{1}{2}} \cdot \Y_{f, t}
\end{align}
where $\phinn_{f}^{\frac{1}{2}}$ is the lower triangular matrix from the Cholesky decomposition.

Speaker \gls{PSD} matrix after Spatial Whitening:
\begin{align}
{\vecG{\Phi}_{\tilde{X}\tilde{X}}}_f &= \frac{\sumTc\sumD M^{(X)}_{f, \check{t}, d}\Ytilde_{f, \check{t}} \cdot \Ytilde_{f, \check{t}}\hermite}{\sumTc \sumD M^{(X)}_{f, \check{t}, d}}.
\end{align}	

PCA beamformer after spatial whitening:
\begin{align}
	\w^{(\mathrm{SW})}_{f} &= \principalComponent{\phixxtilde_f}.
\end{align}	

GEV beamformer:
\begin{align}
\w^{(\mathrm{GEV})}_{f} &= \frac{{\vecG{\Phi}_{NN}}_f^{-\frac{1}{2}\mathrm{H}} \w^{(\mathrm{SW})}_{f}}{\left|\left|{\vecG{\Phi}_{NN}}_f^{-\frac{1}{2}\mathrm{H}} \w^{(\mathrm{SW})}_{f}\right|\right|}. \label{eq:appendix:gev_sw}
\end{align}	

\begin{proof}
	\newcommand{\wGEVf}{\w^{(\mathrm{GEV})}_{f}}
	\newcommand{\wSWf}{\w^{(\mathrm{SW})}_{f}}
	
	GEV criterion:
	\begin{align}
		\wGEVf &= \argmax{\w_f}{\frac{\w_f\hermite\phixx_f\w_f}{\w_f\hermite\phinn_f\w_f}}
	\end{align}
	results in the generalized eigenvalue problem
	\begin{align}
		\phixx_f\wGEVf &= \lambda_\mathrm{max} \phinn_f\wGEVf.
	\end{align}
	Left multiplying both sides with $\phinn_f\inv = \left(\phinn_f^{\frac{1}{2}}\phinn_f^{\frac{1}{2}\mathrm{H}}\right)\inv = \phinn_f^{-\frac{1}{2}\mathrm{H}} \phinn_f^{-\frac{1}{2}}$ and insert \eqref{eq:appendix:gev_sw}
	\begin{align}
		\phinn_f^{-\frac{1}{2}\mathrm{H}} \phinn_f^{-\frac{1}{2}}\phixx_f{\vecG{\Phi}_{NN}}_f^{-\frac{1}{2}\mathrm{H}} \wSWf &= \lambda_\mathrm{max} {\vecG{\Phi}_{NN}}_f^{-\frac{1}{2}\mathrm{H}} \wSWf, \\
		\underbrace{\phinn_f^{-\frac{1}{2}}\phixx_f{\vecG{\Phi}_{NN}}_f^{-\frac{1}{2}\mathrm{H}}}_{\stackrel{!}{=}\phixxtilde_f} \wSWf &= \lambda_\mathrm{max} \wSWf,
	\end{align}
	
	Using the definition for ${\vecG{\Phi}_{\tilde{X}\tilde{X}}}_f$ 
	
	\begin{align}
		{\vecG{\Phi}_{\tilde{X}\tilde{X}}}_f &= \frac{\sumTc\sumD M^{(X)}_{f, \check{t}, d}\phinn_{f}^{-\frac{1}{2}} \Y_{f, t} \Y_{f, t}\hermite \phinn_{f}^{-\frac{1}{2}\mathrm{H}} }{\sumTc \sumD M^{(X)}_{f, \check{t}, d}} \\
		&= \phinn_{f}^{-\frac{1}{2}} \phixx \phinn_{f}^{-\frac{1}{2}\mathrm{H}}
	\end{align}	
	shows, that both are equal.
	
\end{proof}

    \backmatter
    \titleformat{\chapter}
    [display] 
    {\bfseries\huge} 
    {} 
    {0ex}
    {\vspace{-5ex}\titlerule\vspace{1.5ex}\filright~}
    [\vspace{1ex}\titlerule]
    

\chapter{Symbols}
\label{chapter:symbols}

\begin{flushleft}
	\rarray{1.3}

  \begin{longtable}{L{4.2cm}L{\textwidth-5.2cm}} 
   $(\cdot)\transposed$ & transposed \\
   $(\cdot)\transposedinv$ & transposed and invert (swappable)\\
   $(\cdot)\hermite$ & conjugate transposed \\
   $(\cdot)\hermiteinv$ & conjugate transposed and invert (swappable)\\
   $|\cdot|$ & absolute value of a scalar\\
   $||\cdot||$ & Euclidean norm of a vector\\
   $\trace{\cdot}$ & Trace-Operator \\
   $ \cdot \hadamard \cdot$ & Hadamard product or elementwise product \\
   $\Re\{\cdot\}$ & Real Part \\
   $\Im\{\cdot\}$ & Imaginary Part \\
   $\mathrm j$ & Complex number\\
   $\mathrm e$ & Euler's number \\
   $\delta_{nm}$ & Kronecker delta \\
   $\vec{I}$ & Identity matrix \\
   $\vec{0}$ & vector/matrix filled with zeros \\
   $\vec{1}$ & vector/matrix filled with ones \\
   $\lltriangle[very thick]$ & lower triangular matrix  filled with ones \\
   $\urtriangle[very thick]$ &  upper triangular matrix  filled with ones \\
   $\forall_{\text{case}}$ & for all where case is true \\
   $(\cdot)\in\{\dots\}$ & $(\cdot)$ is in the set $\{\dots\}$ \\
   $(\cdot)\in\llblacktriangle$ & a lower triangular matrix \\
   $(\cdot)\in\urblacktriangle$ & a upper triangular matrix \\
   \\
   $d$ & channel index \\
   $n$ & time index \\
   $t$, $\check{t}$ & time index of \gls{stft} \\
   $f$ & frequency index \\
   $\vec{X}_{f,t}$, $\vec{X}_{t}$ & multichannel clean speech in \gls{stft} domain \\
   $X_{f,t,1}$ & reference speech signal on channel one \\
   $\hat{X}_{f,t,1}$ & estimated reference speech signal on channel one \\
   $\vec{N}_{f,t}$, $\vec{N}_{t}$ & multichannel noise in \gls{stft} domain \\
   $\vec{Y}_{f,t}$, $\vec{Y}_{t}$ & multichannel noisy signal in \gls{stft} domain \\
   $M^{(X)}_{f,t,d}$, $M^{(X)}_{t}$ & speech mask (interpretable as speech presence probability) \\
   $M^{(N)}_{f,t,d}$, $M^{(N)}_{t}$ & noise mask (interpretable as noise presence probability) \\
   ${\vecG{\Phi}_{XX}}_f$ & speech \gls{PSD} \\
   ${\vecG{\Phi}_{NN}}_f$ & noise \gls{PSD} \\
   $\vecG{\Phi}$, $\vecG{\Phi}_f$ & noise \gls{PSD}  inverse multiplied with speech \gls{PSD} \\
   $P^{(X)}$ & power of (beamformered) speech signal \\
   $P^{(N)}$ & power of (beamformered) noise signal \\ 
   $\vec{w}_f$, $\vec{w}$ & beamforming vector \\
   $\vec{w}_f^{(\MVDR)}, \vec{w}_f^{(\PCA)}, \vec{w}_f^{(\GEV)}$ & \MVDR/\PCA/\GEV beamforming vector \\
   $\vec{W}$, $\vec{W}_f$ & Eigenvectors \\
   $\vecG{\Lambda}, \vecG{\Lambda}_f$ & Eigenvalue matrix \\
   $\lambda$ & Eigenvalue \\
   $\principalComponent{\cdot}$ & principal component (eigenvector corresponding to the largest eigenvalue) \\
   \\
   $z$ & scalar, input value of a chain \\
   $\zeta$ & placeholder for $z$ and $z\conj$ \\
   $\vec{z}$ & vector, input value of a chain \\
   $\vec{Z}$ & matrix, input value of a chain \\
   $x$ & real part of $z$ \\
   $y$ & imaginary part of $z$ \\
   $s$ & intermediate value of a chain \\
   $\sigma$ & real part of $s$ \\
   $\omega$ & imaginary part of $s$ \\
   $J$ & final real value of a chain, objective function \\
   $f$, $g$ & function \\
   $A$, $B$ & matrix, previous intermediate value of a chain \\
   $C$ & matrix, later intermediate value of a chain \\
   $\der{J}{x}$, $\der{J}{y}$ & gradient, real scalar, partial derivative of $J$ \wrt $x$ or $y$ \\
   $\der{\sigma}{x}$ & real scalar, partial derivative of $\sigma$ \wrt $x$ \\
   $\der{J}{z\conj}$, $\nabla_{z\conj}$ & gradient, scalar, partial derivative of $J$ \wrt $z\conj$ \\
   $\der{J}{z}$, $\nabla_{z}$ & gradient, scalar, partial derivative of $J$ \wrt $z$ \\
   $\der{J}{\vec{z}\conj}$, $\nabla_{\vec{z}\conj}$ & gradient, column vector, partial derivative of $J$ \wrt $\vec{z}\conj$, shape like $\vec{z}$ \\
   $\der{J}{\vec{Z}\conj}$, $\nabla_{\vec{Z}\conj}$ & gradient, matrix, partial derivative of $J$ \wrt $\vec{z}\conj$, shape like $\vec{Z}$ \\
   $\der{\vec{C}}{z\conj}$ & gradient, matrix, partial derivative of $\vec{C}$ \wrt $z\conj$, shape like $\vec{C}$ \\
   $\vec{E}$ & matrix of eigenvalue differences $E_{i,j} = \lambda_j - \lambda_i$\\
   $\vec{F}$ & matrix of inverse eigenvalue differences $F_{i,j} = \left(\lambda_j - \lambda_i\right)$, except on the diagonal where it is zero\\
   $\dft{\cdot}$ & \glsdesc{dft} from time index $n$ to frequency index $f$\\
   $\idft{\cdot}$ & \glsdesc{idft} from frequency index $f$ to time index $n$ \\
   $\rdft{\cdot}$ & \glsdesc{rdft} from time index $n$ to frequency index $f$ \\
   $\irdft{\cdot}$ & \glsdesc{irdft} from frequency index $f$ to time index $n$ \\
  \end{longtable}
\end{flushleft}

    \listoffigures
    \listoftables
    
    \printglossaries
    
    \printbibliography
\end{document}